\documentclass[11pt,a4j,twoside]{article}

\usepackage{amsmath,amssymb}
\usepackage{amsthm}
\usepackage[abbrev]{amsrefs}
\usepackage{latexsym}
\usepackage{graphicx}

\usepackage[english]{babel}
\usepackage{fancyhdr}
\usepackage{a4wide}

\usepackage{bm}

\usepackage[pagewise, mathlines]{lineno}

\allowdisplaybreaks[1]

\numberwithin{equation}{section}

\newtheorem{thm}{Theorem}[section]

\newtheorem{lem}[thm]{Lemma}
\newtheorem{cor}[thm]{Corollary}
\newtheorem{dfn}[thm]{Definition}
\theoremstyle{definition} 

\newtheorem{rem}[thm]{Remark}

\newcommand\ND{\newcommand}

\ND\lref[1]{Lemma~\ref{#1}}
\ND\tref[1]{Theorem~\ref{#1}}
\ND\pref[1]{Proposition~\ref{#1}}
\ND\sref[1]{Section~\ref{#1}}
\ND\ssref[1]{Subsection~\ref{#1}}
\ND\aref[1]{Appendix~\ref{#1}}
\ND\rref[1]{Remark~\ref{#1}}
\ND\cref[1]{Corollary~\ref{#1}}
\ND\eref[1]{Example~\ref{#1}}
\ND\fref[1]{Fig.\ {#1} }
\ND\lsref[1]{Lemmas~\ref{#1}}
\ND\tsref[1]{Theorems~\ref{#1}}
\ND\dref[1]{Definition~\ref{#1}}
\ND\psref[1]{Propositions~\ref{#1}}
\ND\rsref[1]{Remarks~\ref{#1}}
\ND\sssref[1]{Subsections~\ref{#1}}
\ND\esref[1]{Examples~\ref{#1}}
\ND\asref[1]{Assumption~\ref{#1}}

\newcommand{\ep}{\varepsilon}

\newcommand{\wilde}{\widetilde}

\newcommand{\h}{\quad}
\newcommand{\dis}{\displaystyle}

\newcommand{\cl}{c\`adl\`ag\ }
\newcommand{\Prob}{\mathbb{P}}
\newcommand{\Ex}{\mathbb{E}}
\newcommand{\bi}{\bar{\iota}}

\newcommand{\bp}{\mathbf{p}}

\newcommand{\bZ}{\mathbf{Z}}
\newcommand{\bdm}{\mathbf{m}}
\newcommand{\cX}{\mathcal{X}}

\newcommand{\cF}{\mathcal{F}}
\newcommand{\cE}{\mathcal{E}}
\newcommand{\cA}{\mathcal{A}}

\newcommand{\cL}{\mathcal{L}}
\newcommand{\cN}{\mathcal{N}}

\newcommand{\cP}{\mathcal{P}}

\newcommand{\bbR}{\mathbb{R}}
\newcommand{\sS}{\mathsf{S}}

\title{On sequences of martingales with jumps on Riemannian submanifolds}
\date{}
\author{Fumiya Okazaki}
\begin{document}
\maketitle
\footnote{Email address: okazaki.f.660b@m.isct.ac.jp}
\renewcommand{\thepage}{\arabic{page}}
\begin{abstract}
In this article, we investigate sequences of discontinuous martingales on submanifolds of higher-dimensional Euclidean space. Those sequences naturally arise when we deal with a sequence of harmonic maps with respect to non-local Dirichlet forms, such as fractional harmonic maps. We prove that the semimartingale topology is equivalent to the topology of locally uniform convergence in probability on the space of discontinuous martingales on manifolds with bounded jumps. In particular, we show that the limit of any sequence of discontinuous martingales on a compact Riemannian manifold, with respect to the topology of locally uniform convergence in probability, is a martingale on the manifold. As an application, we prove that for sequences of quasi-harmonic maps on an open set with respect to a non-local Dirichlet form, local $L^{\infty}$-convergence in the open set with pointwise convergence q.e. on the whole space implies strong local convergence in the Dirichlet space and the limit preserves harmonicity.
\end{abstract}

\begin{flushleft}
{\bf MSC2020 Subject Classifications: 60G44, 60H05}
\end{flushleft}

\section{Introduction}
Martingales on manifolds have been studied in relation to harmonic maps which are defined as critical points of the Dirichlet energy defined for maps between two Riemannian manifolds. In fact, it is known that the image of Brownian motion on the domain Riemannian manifold by a harmonic map is a continuous martingale on the target manifold. Here the notion of martingales on manifolds is defined by continuous semimartingales along which the It\^o integral of each 1-form is a local martingale.

Recently, the notion of harmonic maps with respect to the fractional Laplacian, which are called fractional harmonic maps, was introduced in \cites{DaLioRiv11, DaLioRiv112} and their regularity has been studied in \cites{DaLioRiv11, DaLioRiv112, MazSch18, MP20, MSire15, MPS21}. More generally, harmonic maps with respect to non-local Dirichlet forms also admit a probabilistic characterization through stochastic processes. In fact, it has been shown in \cite{Oka24} that a map valued in a submanifold of Euclidean space solves the Euler-Lagrange equation for a Dirichlet form if and only if the image of the Markov process associated with the Dirichlet form by the map is a discontinuous martingale on the target submanifold. Here the notion of discontinuous martingales on manifolds is based on \cite{Pic91}. In the same way as the continuous case, discontinuous martingales on manifolds can also be defined through the It\^o integral of 1-forms along semimartingales on manifolds. Here we remark that while the It\^o integral along continuous semimartingales on manifolds can be defined through linear connections on tangent bundles, they are not sufficient ingredients in the discontinuous case. 
In \cite{Pic91}, the notion of connection rules has been introduced to define It\^o integral for discontinuous semimartingales on manifolds.

In this article, we focus on sequences of discontinuous martingales on Riemannian submanifolds of the higher-dimensional Euclidean space. Such sequences naturally appear if we consider sequences of harmonic maps.
In the case of continuous martingales, it is known that the limit of a sequence of continuous martingales on a manifold in the topology of locally uniform convergence in probability (u.c.p.) is also a continuous martingale on the manifold (e.g. Theorem (4.43) of \cite{Emery89}). In the same way as sequences of continuous martingales, the limit of a sequence of discontinuous martingales is expected to be a martingale on the manifold. The aim of this article is to verify such a statement under some assumptions on jumps. However, there are some obstacles that arise from manifolds and jumps. Since discontinuous martingales on manifolds are defined through the It\^o integral, we need to deal with sequences of stochastic integrals. However, in general, the almost sure uniform convergence of a sequence of semimartingales does not yield the convergence in probability of the sequence of the stochastic integrals induced from the sequence of the semimartingales. Therefore, we need to verify the convergence in a stronger topology for manifold-valued martingales in order to show the convergence of the sequence of the stochastic integrals. In \cite{AT98}, the semimartingale topology was studied and the stability of stochastic differential equations was also studied. In the same article, it was shown that if a sequence of continuous martingales on a manifold converges in the u.c.p. topology, then it converges to a continuous martingale on the manifold in the semimartingale topology. 
Another obstacle is the fact that the space of martingales on the Euclidean space is not closed with respect to the semimartingale topology without an additional assumption on the jumps. An example is given in \cite{Emery79}.

Our first result is the equivalence between the two kinds of topologies defined on the space of manifold-valued martingales which may have jumps. One of them is the u.c.p. topology and the other is the topology induced by the so-called $\mathbb{H}^p$-norm. We recall the precise definition of $\mathbb{H}^p$-norm in \sref{Preliminaries}. \tref{smtopology} below shows the Burkholder-Davis-Gundy-type inequality for discontinuous martingales on manifolds. Burkholder-Davis-Gundy's inequality for continuous martingales on manifolds was shown in \cite{Darling96} in some situations. In the study of continuous martingales on manifolds including Burkholder-Davis-Gundy's inequality, the following property of martingales plays an important role: The images of continuous martingales on a small open set by convex functions on the open set are local submartingales.
However, since martingales with jumps possibly exit an open set by large jumps and non-trivial convex functions do not globally exist on manifolds in general, discontinuous martingales do not admit this property in our setting. 
Therefore, we employ another method using the exponential function on $\mathbb{R}^d$ and the scaling of processes to obtain the Burkholder-Davis-Gundy-type inequality in \sref{pfsm}.\\

Let $M$ be a properly embedded complete Riemannian submanifold of the higher dimensional Euclidean space $\mathbb{R}^d$, i.e., there exists a differentiable injective map $\iota \colon M \to \mathbb{R}^d$ such that the map $\iota \colon M \to \iota(M)$ is a diffeomorphism onto its image endowed with the relative topology and the inverse image of every compact set in $\mathbb{R}^d$ is compact. We say that a semimartingale $X$ on $M$ is a martingale on $M$ if the It\^o integral of each smooth vector field $\mathcal{X}$ on $M$ along $X$ is a local martingale, where we consider the It\^o integral by regarding both processes $X$ and $\mathcal{X}(X)$ as $\mathbb{R}^d$-valued processes. This definition coincides with the notion of $\gamma$-martingales introduced in \cite{Pic91} for the connection rule $\gamma \colon M \times M \to TM$,
\begin{align}\label{ConneEmbed}
\gamma (x,y) = \Pi_x(y-x),
\end{align}
where $\Pi_x \colon \mathbb{R}^d \to T_xM$ is the orthonormal projection. In order to incorporate the inside killing of martingales, we employ a slightly extended definition of martingales as considered in \cites{Oka23, Oka24}. We fix a point $\bp \in \mathbb{R}^d \backslash M$ and denote $M$-valued martingales with end points by triplets $(X,\zeta,\bp)$, where $\zeta$ is a stopping time which represents the killing time of the martingale $X$ and $\bp$ represents the trap. The killing time is essential when we consider martingales associated with harmonic maps with respect to Dirichlet forms with a killing term. We will provide the precise definition of martingales on $M$ in Definitions \ref{KilledSemi} and  \ref{ExtendedInt}. For $f\in C^{\infty}(M)$, we denote by $\bar{f}$ an extension of $f$ to a function in $C^{\infty}(\mathbb{R}^d)$ satisfying
\begin{align}\label{nabla}
\nabla f(x)=D\bar{f}(x)\ \text{on}\ M,
\end{align}
where $D\bar{f}$ stands for the gradient on $\mathbb{R}^d$ (such an extension exists since the embedding is proper). We further introduce the extension of the embedding $\iota$ denoted by $\bi = (\bi^1,\dots,\bi^d)\in C^{\infty}(\mathbb{R}^d ; \mathbb{R}^d)$ such that $\bi$ satisfies $\bi(\bp)=\bp$ in addition to \eqref{nabla}. For a point $\bp \in \mathbb{R}^d \backslash M$ and $\alpha, \beta >0$, we denote by $\mathcal{M}_{\alpha,\beta}(M,\bp)$ the set of $M$-valued martingales $(X,\zeta,\bp)$ with the end point $\bp$ satisfying
\begin{align}
\sup_{0\leq s <\infty} |\Delta \bi(X_s)| \leq \alpha, \label{jumpbound}\\
|X_0|\leq \beta \label{initialbound}.
\end{align}
For a stochastic process $H$ and a stopping time $\tau$, we write the stopped process as $H^{\tau}$ defined by
\[
H^{\tau}_t(\omega)=H_{t\land \tau (\omega)}(\omega).
\]
\begin{thm}\label{smtopology}
Let $M$ be a properly embedded complete Riemannian submanifold of the higher dimensional Euclidean space $\mathbb{R}^d$ and $\bp$ a point in $\mathbb{R}^d \backslash M$. Let $\alpha,\beta >0$. Let $(X,\zeta,\bp)$ and $(Y,\zeta',\bp)$ be $M$-valued martingales with the end point $\bp$ in $\mathcal{M}_{\alpha,\beta}(M,\bp)$.
Then there exists $R=R(M,\alpha,\beta)>0$ such that for all $p\in [1,\infty)$, there exists $C=C(M,R,p,\alpha,\beta)>0$ such that 
\[
\| \bi(X)^{t\land \tau}-\bi(Y)^{t\land \tau} \|_{\mathbb{H}^p} \leq C \Ex \left[ \left( \sup_{0\leq s \leq t\land \tau}|\bi(X_s)-\bi(Y_s)| \right)^{2p} \right]^{\frac{1}{4p}}
\]
for every $t\geq 0$ and stopping time $\tau$ with
\begin{align}
\sup_{0\leq s < \tau} |\bi(X_s)-\bi(X_0)| &\leq R,\label{smallrange}\\
\sup_{0\leq s < \tau} |\bi(Y_s)-\bi(Y_0)| &\leq R.\label{smallrangeY}
\end{align}
\end{thm}
We note that if $M$ is compact, conditions \eqref{jumpbound} and \eqref{initialbound} are automatically satisfied.
Since the $\mathbb{H}^p$-norm is related to the semimartingale topology, we can show the equivalence of the u.c.p. topology and the semimartingale topology on the space of manifold-valued martingales as a corollary of  \tref{smtopology}.
\begin{cor}\label{equitopo}
Let $M$ be a second-countable properly embedded complete Riemannian submanifold of the higher dimensional Euclidean space $\mathbb{R}^d$. We fix $\bp \in \mathbb{R}^d \backslash M$ and $\alpha, \beta>0$. Let $(X^n,\zeta^n,\bp)$ be a sequence in $\mathcal{M}_{\alpha,\beta}(M, \bp)$.
Assume that $\{X^n\}_{n\in \mathbb{N}}$ converges to a process $X$ on $\mathbb{R}^d$ with respect to the u.c.p. topology.
Then $X$ is a semimartingale and the sequence $\{X^{n}\}_{n\in \mathbb{N}}$ converges to $X$ with respect to the semimartingale topology.
Moreover, there exists a subsequence $\{n_k\}_{k\in \mathbb{N}}$ and a stopping time $\zeta$ such that
\begin{align}\label{liminfzeta}
\Prob \left( \bigcap_{l=1}^{\infty} \liminf_{k\to \infty} \{ \zeta^{n_k} \land l=\zeta \land l\} \right)=1
\end{align}
and the triplet $(X,\zeta,\bp)$ is a martingale on $M$ with the end point $\bp$.
\end{cor}
The event in \eqref{liminfzeta} means that for almost all $\omega \in \{ \zeta < \infty\}$, there exists $N(\omega) \in \mathbb{N}$ such that $\zeta^{n_k}(\omega)=\zeta(\omega)$ for all $k\geq N(\omega)$.
\begin{rem}
The assumption \eqref{jumpbound} is essential since as mentioned above, an example in \cite{Emery79} shows that the limit of a sequence of local martingales on $\mathbb{R}$ with respect to the semimartingale topology is not necessarily a local martingale on $\mathbb{R}$.
\end{rem}

As an application of \cref{equitopo}, we obtain a strong local $\cE_1$-convergence result for sequences of harmonic maps associated with non-local Dirichlet forms under the local $L^{\infty}$-convergence. For harmonic maps with respect to local Dirichlet forms, this property can be established by the characterization by local convex functions and the resulting subharmonicity of their compositions with the map. However, in the present setting of non-local Dirichlet forms, such a localization is no longer available directly. The significance of \cref{equitopo} in this application is that we can make use of the probabilistic characterization, which is still valid in non-local cases with various target manifolds.

We briefly describe the setting for harmonic maps. For full details, see \sref{PfofSeqharm} or refer to \cites{FOT, ChenFuku}. Let $E$ be a locally compact separable metric space, $\bdm$ a positive Radon measure with full support on $E$. We add a cemetery point $\partial$ to the space $E$ and extend any function $u \colon E \to \mathbb{R}$ to $E \cup \partial$ by the convention $u(\partial)=0$. Let $(\cE,\cF)$ be a regular Dirichlet form on $L^2(E;\bdm)$ and let $\bZ=(\Omega, \{Z_t\}_{t\geq 0}, \zeta, \{\Prob_z\}_{z\in E})$ be an $\bdm$-symmetric Hunt process on $E$. Let $\{ \cF^Z_t \}_{t\geq 0}$ be the minimal admissible augmented filtration generated by $Z$. For an open set $D$, we set the first exit time
\[
\tau_D:=\inf \{ t> 0 \mid Z_t \notin D \}.
\]
Let $(\cE^D,\cF^D)$ be the part form of $(\cE,\cF)$ on $D$.
\begin{dfn}
Assume $0 \notin M$.
Let $\bZ=(\{Z_t\}_{t\geq 0}, \zeta, \{\Prob_z\}_{z\in E})$ be an $\bdm$-symmetric Hunt process on $E$.
Let $D \subset E$ be an open set.
Let $u\colon E \to \mathbb{R}^d$ be a Borel measurable map. We say that $u$ is a quasi-harmonic map on $D$ with respect to $\bZ$ if for each $D_1 \Subset D$, $(\{u(Z_{t\land \tau_{D_1}})\}_{t\geq 0},\hat{\zeta},0)$ is a $(\Prob_z, \{\cF^Z_t\}_{t\geq 0})$-martingale on $M$ with an end point for q.e. $z\in E$, where
\[
\hat{\zeta}=\zeta \mathbf{1}_{\{ \tau_{D_1} \geq \zeta \}} +\infty \mathbf{1}_{\{ \tau_{D_1}< \zeta \}}.
\]
\end{dfn}
According to Theorem 1.1 of \cite{Oka24}, the notion of quasi-harmonic map is equivalent to weakly harmonic map, which is defined by the Euler-Lagrange equation derived by the variation principle as long as $u \in \cF^D_{loc}(M)$, where $\cF^D_{loc}(M)$ is the $M$-valued local Dirichlet space.
\begin{thm}\label{harmonicseq}
Let $D \subset E$ be an open set.
Let $M$ be a properly embedded compact Riemannian submanifold of $\mathbb{R}^d$ with $0 \notin M$.
Let $u_n \in \cF^D_{loc}(M)$ be a quasi-harmonic map on $D$ for $n \in \mathbb{N}$. Assume that there exists a map $u\colon E \to \mathbb{R}^d$ such that
\[
u_n \to u\ \text{in}\ L^{\infty}_{loc}(D)\ \text{and pointwise q.e. on}\ E.
\]
Then $u \in \cF^D_{loc}(M)$, $u_n \to u$ strongly in $\cF^D_{loc}$, and $u$ is also a quasi-harmonic map on $D$. Here, the strong convergence in $\cF^D_{loc}$ is defined by
\[
\cE_1^D(\phi (u_n-u), \phi (u_n-u)) \to 0
\]
as $n \to \infty$ for any $\phi \in C_0(D) \cap \cF$. Here, for $v\in \cF(\mathbb{R}^d)$, we set
\[
\cE(v,v):=\sum_{i=1}^d\cE(v^i,v^i).
\]
\end{thm}
We give a comment on \tref{harmonicseq}. For sphere-valued stationary or energy-minimizing harmonic maps with respect to a fractional Laplacian, the compactness in local fractional Sobolev spaces under the uniform energy bound has been studied in \cites{MPS21}. In contrast, in \tref{harmonicseq} above, the uniform bound for the energy is not required and we claim that the strong local $\cE_1$-convergence is obtained, although it assumes local $L^{\infty}$-convergence in $D$. Moreover, the assumption that $u \in \cF^D_{loc}(M)$ is not required either.\\

We give an outline of the paper. In \sref{Preliminaries}, we first recall the notion of semimartingale topology and some facts regarding it. Mainly we refer to \cites{Protter05} there. We also recall the It\^o calculus for discontinuous semimartingales on manifolds. In \sref{pfsm}, we provide the proofs of \tref{smtopology} and \cref{equitopo}. In \sref{PfofSeqharm}, we recall some notions regarding Dirichlet forms and give a proof of \tref{harmonicseq}.\\

Throughout this paper, we fix a filtered probability space $(\Omega, \mathcal{F}, \{\mathcal{F}_t\}_{t\geq 0}, \Prob)$ satisfying the usual conditions.
We denote the set of $\{\mathcal{F}_t \}$-predictable processes by $\mathcal{P}$. We denote the set of semimartingales on $\mathbb{R}^d$ by $\mathcal{S}(\mathbb{R}^d)$. For a \cl process $X$, we denote by $X_-$ the process obtained by taking the left-limit of $X$, namely,
\[
X_{t-}(\omega)=\lim_{s\nearrow t}X_s(\omega).
\] 
\section{Preliminaries}\label{Preliminaries}
In this section, we recall definitions of some norms on $\mathcal{S}(\mathbb{R}^d)$ and their properties. For a \cl process $X$, we set
\[
r(X):=\sum_{k=1}^{\infty} \frac{1}{2^k} \Ex \left[ 1\land \sup_{0\leq t \leq k}|X_t| \right].
\]
Then for a sequence $X^n$ of \cl processes, $r(X^n)\to 0$ as $n \to \infty$ if and only if $X^n$ converges to $0$ in u.c.p. Next we suppose $X \in \mathcal{S}(\mathbb{R}^d)$. We set
\[
\hat r(X):=\sup_{H\in \mathcal{P}_1^d}r\left( \int \langle H, dX \rangle \right),
\]
where
\begin{align*}
\mathcal{P}_1&:=\{ H\in \mathcal{P} \mid |H_t|\leq 1\ \text{for all}\ t\geq 0\ \text{a.s.}\},
\end{align*}
and $\mathcal{P}_1^d$ is the Cartesian product of $d$-copies of $\mathcal{P}_1$. Then $\hat r$ induces a topology on $\mathcal{S}(\mathbb{R}^d)$ and the topology on $\mathcal{S}(\mathbb{R}^d)$ is called the semimartingale topology.

Next we recall another norm defined for semimartingales and its relation to $\hat r$. Let $p\in [1,\infty)$. 
For an $\mathbb{R}$-valued semimartingale $X$, we set
\begin{align*}
\| X \|_{\mathbb{H}^p}&:=\inf_{X=X_0+M+A} \Ex \left[ |X_0|^p+[M,M]_{\infty}^{\frac{p}{2}}+\left(\int_0^{\infty}|dA_t|\right)^p  \right]^{\frac{1}{p}},\\
\mathbb{H}^p&:=\{ X\in \mathcal{S}(\mathbb{R}) \mid \|X\|_{\mathbb{H}^p}<\infty \},
\end{align*}
where the infimum is taken over all decompositions of the semimartingale $X$ into a local martingale $M$ and a process $A$ of locally finite variation with $M_0=A_0=0$. For an $\mathbb{R}^d$-valued semimartingale $X$, we set
\[
\| X \|_{\mathbb{H}^p}:=\left( \sum_{k=1}^d \| X^k \|_{\mathbb{H}^p}^p \right)^{\frac{1}{p}}.
\]
\begin{rem}\label{Hpconv}
In Theorem 2 of \cite{Emery79} (cf: Theorem 14 in p. 264 of \cite{Protter05}), the relationship between the semimartingale topology and the norm $\| \cdot \|_{\mathbb{H}^p}$ is clarified. Let $p\in [1,\infty)$, $X^n, X\in \mathcal{S}(\mathbb{R})$. Then the following hold.
\begin{itemize}
\item[(1)]If $X^n$ converges to $X$ in the semimartingale topology, then there exists a subsequence $X^{n_k}$ and a non-decreasing sequence of stopping times $\sigma_m$ which tends to $\infty$ such that $X^{n_k,\sigma_m-}$ converges to $X^{\sigma_m-}$ in $\mathbb{H}^p$ as $k \to \infty$ for each $m$.
\item[(2)]If there exists a non-decreasing sequence of stopping times $\sigma_m$ which tends to $\infty$ such that $X^{n,\sigma_m-}$ converges to $X^{\sigma_m-}$ in $\mathbb{H}^p$, then $X^n$ converges to $X$ in the semimartingale topology.
\end{itemize}
\end{rem}

Next we recall the It\^o calculus on manifolds. Let $M$ be a manifold. The definition of semimartingales on $M$ is simple as follows.
\begin{dfn}
Let $X$ be an $M$-valued \cl process. The process $X$ is called an $M$-valued semimartingale if for all $f\in C^{\infty}(M)$, $f(X)$ is a semimartingale on $\mathbb{R}$.
\end{dfn}
In order to define the It\^o integral, we need to determine jumps of semimartingales as tangent vectors. 
\begin{dfn}
Let $\gamma:M\times M\to TM$ be a measurable map and suppose that $\gamma$ is $C^{\infty}$ on a neighborhood of $\mathrm{diag}(M)$. Then $\gamma$ is called a connection rule if it satisfies the following conditions, for all $x,y\in M$,
\begin{enumerate}
\item[(i)] $\gamma (x,y)\in T_xM$;
\item[(ii)] $\gamma (x,x)=0$;
\item[(iii)] $(d \gamma (x,\cdot))_x=id_{T_xM}.$
\end{enumerate}
\end{dfn}
For a given semimartingale $X$ on $M$ and a connection rule $\gamma$, we can regard jumps of $X$ as $\gamma(X_{t-},X_t) \in T_{X_{t-}}M$. Then for each smooth 2-tensor field $b$ on $M$, we can define the quadratic variation $\dis \int b(X_-)\, \gamma d[X,X]$ (see \cite{Pic91} or \cite{Oka23} for details). Denote the continuous locally bounded variation part by $\dis \int b(X_-)\, d[X,X]^c$, which is independent of the choice of the connection rule. We can also define the It\^o integral of each 1-form $\alpha$ on $M$ denoted by $\dis \int \alpha(X_-)\, \gamma dX$ in such a way that it holds that
\begin{align*}
\int df(X_-)\, \gamma dX&=f(X)-f(X_0)-\frac{1}{2}\int \nabla df (X_-)\, d[X,X]^c\\
&\h -\sum_{0<s\leq \cdot} \{f(X_s)-f(X_{s-})-\langle df(X_{s-}),\gamma(X_{s-},X_s) \rangle \}
\end{align*}
for every $f\in C^2(M)$.
In this article, we focus on the case where the manifold $M$ is a Riemannian submanifold of the higher dimensional Euclidean space and the connection rule $\gamma$ is given by \eqref{ConneEmbed}. In \cite{Oka24}, in order to incorporate the inside killing of semimartingales, the author used the following slightly extended definition of semimartingales on manifolds.
\begin{dfn}\label{KilledSemi}
Let $M$ be an isometrically and properly embedded complete Riemannian submanifold of $\mathbb{R}^d$. Let $X$ be an $\mathbb{R}^d$-valued semimartingale, $\zeta$ a stopping time, and $\bp$ a point in $\mathbb{R}^d \backslash M$. The triplet $(X,\zeta,\bp)$ is called an $M$-valued semimartingale with an end point if $X_t \in M$ for $t\in [0,\zeta)$ and $X_t=\bp$ for $t\geq \zeta$.
\end{dfn}
\begin{rem}
Since $M$-valued semimartingales with end points are semimartingales on $\mathbb{R}^d$, $X_{\zeta-}$ exists on $\{ \zeta < \infty \}$. Moreover, since we assume that $M$ is properly embedded in $\mathbb{R}^d$, $M$ is closed in $\mathbb{R}^d$. This yields $X_{\zeta-} \in M$ on $\{ \zeta < \infty \}$. In particular, the finite-time explosion is excluded in the present framework.
\end{rem}

\begin{dfn}\label{ExtendedInt}
Let $M$ be an isometrically and properly embedded complete Riemannian submanifold of $\mathbb{R}^d$ and $(X,\zeta,\bp)$ an $M$-valued semimartingale with an end point. Then $(X,\zeta,\bp)$ is called a martingale with an end point if for each $\mathcal{X}\in \mathfrak{X}(M)$, the stochastic integral $\dis \int \langle \mathcal{X}(X_-),dX\rangle$ on $\mathbb{R}^d$ is a local martingale, where the set $\mathfrak{X}(M)$ stands for the set of all smooth vector fields on $M$.
\end{dfn}
Note that the integral in \dref{ExtendedInt} is well-defined even though $\mathcal{X}$ is not defined at the end point of $X$ since $X$ stops after reaching the end point. Moreover, we can easily check that if $\zeta=\infty$ a.s., then the notion of martingales on $M$ with an end point coincides with the one of $\gamma$-martingales. In fact, for every $f\in C^{\infty}(M)$, it holds that
\begin{align*}
\int df(X_-)\, \gamma dX=\int \langle D\bar{f}(X_-),dX\rangle,
\end{align*}
where $\bar{f}$ is an extension of $f$ satisfying \eqref{nabla}.
\begin{rem}
Typically the Riemannian quadratic variation of a continuous semimartingale $X$ on $M$ is defined by
\[
\int g(X)\, d[X,X],
\]
where $g$ is a Riemannian metric and it coincides with the quadratic variation $[\iota(X),\iota(X)]$ of the $\mathbb{R}^d$-valued semimartingale $\iota(X)$. However, in the case where $X$ has jumps, we need to distinguish these two notions since they do not coincide in general. In \sref{pfsm}, we focus on the quadratic variation $[\bi(X),\bi(X)]$ of a martingale $X$ with an end point in the same way as \cite{Oka23}.
\end{rem}
\section{Proof of \tref{smtopology}}\label{pfsm}
In order to prove \tref{smtopology}, we will prepare some estimates for the quadratic variation of martingales. The following lemma provides the uniform bound of the expectation of quadratic variations of martingales on a small open set. A similar estimate was obtained in \cite{Pic94} for martingales defined through barycenters but we do not assume any geometric condition on $M$ in our setting regarding barycenters. We employ a similar technique to the proof of Theorem 1.2 of \cite{Oka23}, where the composition of scaled martingales and the exponential function was considered. The technique of using the exponential function was employed in \cite{HeYanZheng83} in the study of the convergence of continuous martingales on manifolds and in \cite{Oka23}, the method was modified for discontinuous cases.

Prior to the proofs of lemmas below, we introduce some notation. For $R>0$, we set
\begin{align*}
a_1(R)&:=\sup_{x\in \overline{B_R}}|D\bi (x)|,\\
a_2(R)&:=\sup_{x\in \overline{B_R}}\sup_{u\in \mathbb{R}^d, |u|=1}|\mathrm{Hess}\, \bi(x)(u,u)|,\\
a_3(R)&:=\sup_{x \in \overline{B_R}} \sup_{i,j,k,l=1,\dots,d}|\partial_i \partial_j \partial_k \bi^l (x)|,
\end{align*}
where $B_R$ is the ball in $\mathbb{R}^d$ with the radius $R>0$ centered at $0$.
\begin{lem}\label{Exquad}
Let $\alpha,\beta>0$. Then there exists $R=R(M,\alpha,\beta)>0$ such that for all $p\in [1,\infty)$, there exists $C=C(M,R,p,\alpha,\beta)$ such that for every $(X,\zeta,\bp) \in \mathcal{M}_{\alpha,\beta}(M,\bp)$ and stopping time $\tau$ satisfying $\eqref{smallrange}$ for $R$, it holds that
\[
\Ex \left[ [\bi(X),\bi(X)]^p_{\tau} \right]\leq C.
\]
\end{lem}
\begin{proof}
Let $\sigma$ be an arbitrary stopping time and set $\tau'=\tau \land \sigma$ (We will take a concrete stopping time later in \eqref{sigman}). For an arbitrary $R>0$, we set
\begin{align*}
Y_t&:=\frac{1}{R}\left(\bi(X_t)^{\tau'}-\bi(X_0)\right)\\
Z^i_t&:=e^{Y^i_t},\\
Z_t&=\sum_{i=1}^dZ^i_t.
\end{align*}
Then by It\^o's formula,
\begin{align*}
Z^i_{t}-Z^i_0&=\int_{0+}^te^{Y^i_{s-}}\, dY^i_s +\frac{1}{2}\int_{0+}^te^{Y^i_{s-}}\, d[Y^i,Y^i]^c_s\\
&\h +\sum_{0<s\leq t}\left(e^{Y^i_s}-e^{Y^i_{s-}}-e^{Y^i_{s-}}\Delta Y^i_s \right)\\
&=\frac{1}{R} \int_{0+}^te^{Y^i_{s-}}\, d\bar{\iota^i}(X)^{\tau'}_s +\frac{1}{2R^2}\int_{0+}^te^{Y^i_{s-}}\, d[\bar{\iota^i}(X)^{\tau'},\bar{\iota^i}(X)^{\tau'}]^c_s\\
&\h +\sum_{0<s\leq t}\left(e^{Y^i_s}-e^{Y^i_{s-}}-e^{Y^i_{s-}}\Delta Y^i_s \right)\\
&=\frac{1}{R} \int_{0+}^te^{Y^i_{s-}}\, d\left( \int_{0+}^sD_j \bar{\iota ^i}(X_{u-})^{\tau'}\, d\bar{\iota^j}(X)^{\tau'}_u\right)\\
&\h + \frac{1}{2R}\int_{0+}^t e^{Y^i_{s-}} D_jD_k\bar{\iota^i}(X_{s-})^{\tau'}\, d[\bar{\iota^j}(X)^{\tau'},\bar{\iota^k}(X)^{\tau'}]^c_s\\
&\h +\frac{1}{2R^2}\int_{0+}^te^{Y^i_{s-}}\, d[\bar{\iota^i}(X)^{\tau'},\bar{\iota^i}(X)^{\tau'}]^c_s\\
&\h +\frac{1}{R}\sum_{0<s\leq t}\{ e^{Y^i_{s-}}\left(\bar{\iota^i}(X_{s})^{\tau'}-\bar{\iota^i}(X_{s-})^{\tau'}-\langle D\bar{\iota^i}(X_{s-})^{\tau'},\Delta \bar{\iota} (X_s)^{\tau'} \rangle \right) \} \\
&\h + \sum_{0<s\leq t}\left(e^{Y^i_s}-e^{Y^i_{s-}}-e^{Y^i_{s-}}\Delta Y^i_s \right).
\end{align*}
We set
\begin{align*}
N_t&:=\frac{1}{R} \sum_{i=1}^d\int_{0+}^te^{Y^i_{s-}}D_j\bi^i(X_{s-})^{\tau'}\, d\bi^j(X_s)^{\tau'},\\
J_t&:=\sum_{i=1}^d\frac{1}{R}\mathbf{1}_{\{ \tau' \leq t\}}\left( e^{Y^i_{\tau'-}}\left( \bar{\iota^i}(X_{\tau'}) - \bar{\iota^i}(X_{\tau'-}) - \langle D\bar{\iota^i}(X_{\tau'-}), \Delta \bar{\iota} (X_{\tau'}) \rangle \right)\right),\\
K_t&:=\sum_{i=1}^d\mathbf{1}_{\{ \tau' \leq t\}}\left( e^{Y^i_{\tau'}}-e^{Y^i_{\tau'-}}-e^{Y^i_{\tau'-}}\Delta Y^i_{\tau'} \right),\\
A_t&:=\sum_{i=1}^d\left(\frac{1}{2R}\int_{0+}^t e^{Y^i_{s-}} D_jD_k\bar{\iota^i}(X_{s-})^{\tau'}\, d[\bar{\iota^j}(X_{s-})^{\tau'},\bar{\iota^k}(X_{s-})^{\tau'}]^c_s\right. \nonumber \\
&\h \left. +\frac{1}{2R^2}\int_{0+}^te^{Y^i_{s-}}\, d[\bar{\iota^i}(X)^{\tau'},\bar{\iota^i}(X)^{\tau'}]^c_s \right),\\
B_t&:=\sum_{i=1}^d\left( \frac{1}{R}\sum_{0<s\leq t}\left\{ e^{Y^i_{s-}}\left(\bar{\iota^i}(X_{s})^{\tau'}-\bar{\iota^i}(X_{s-})^{\tau'}-\langle D\bar{\iota^i}(X_{s-})^{\tau'},\Delta \bar{\iota} (X_s)^{\tau'}\rangle \right) \right\} \right. \nonumber \\
&\h \left.+ \sum_{0<s\leq t}\left(e^{Y^i_s}-e^{Y^i_{s-}}-e^{Y^i_{s-}}\Delta Y^i_s \right) \right)-J_t-K_t.
\end{align*}
Then $Z_t-J_t-K_t=N_t+A_t+B_t$ and $N$ is a local martingale since $X$ is an $M$-valued martingale with an end point. We fix $R_0$ with $R\leq R_0$ and set $a_i:=a_i(R_0+\alpha+\beta)$, $i=1,2$. Then the process $A$ satisfies
\begin{align*}
dA_s&=\sum_{i=1}^d \left( \frac{1}{2R}e^{Y^i_{s-}} D_jD_k\bar{\iota^i}(X_{s-})^{\tau'}\, d[\bar{\iota^j}(X)^{\tau'},\bar{\iota^k}(X)^{\tau'}]^c_s+\frac{1}{2R^2}e^{Y^i_{s-}}\, d[\bar{\iota^i}(X)^{\tau'},\bar{\iota^i}(X)^{\tau'}]^c_s \right)\\
&\geq \sum_{i=1}^d\left( -\frac{ea_2}{R}  d[\bar{\iota}(X)^{\tau'}, \bar{\iota}(X)^{\tau'}]^c_s + \frac{e^{-1}}{2R^2} d[\bar{\iota^i}(X)^{\tau'},\bar{\iota^i}(X)^{\tau'}]^c_s\right)\\
&= \left( -\frac{dea_2}{R} + \frac{e^{-1}}{2R^2} \right)d[\bar{\iota}(X)^{\tau'}, \bar{\iota}(X)^{\tau'}]^c_s.
\end{align*}
In a similar way, we have
\begin{align*}
dB_s &= \sum_{i=1}^d\frac{e^{Y^i_{s-}}}{R}\left(\bar{\iota^i}(X_{s})^{\tau'}-\bar{\iota^i}(X_{s-})^{\tau'}-\langle D\bar{\iota^i}(X_{s-})^{\tau'},\Delta \bar{\iota} (X_s)^{\tau'}\rangle \right)\mathbf{1}_{\{s<\tau'\}}\\
&\h +\sum_{i=1}^d\left(e^{Y^i_s}-e^{Y^i_{s-}}-e^{Y^i_{s-}}\Delta Y^i_s\right)\mathbf{1}_{\{s<\tau'\}}\\
&\geq - \frac{de a_2}{R}|\Delta \bar{\iota} (X_s)|^2\mathbf{1}_{\{s<\tau'\}}\\
&\h + \inf \left\{ \frac{e^y-e^x-e^x(y-x)}{(y-x)^2}\mid x,y\in \mathbb{R}, x\neq y,  |x|, |y| \leq 1 \right\} \sum_{i=1}^d|\Delta Y^i_s |^2\mathbf{1}_{\{s<\tau'\}}\\
& \geq \left( -\frac{de a_2}{R} + \frac{e^{-1}}{2R^2}\right) |\Delta \bar{\iota}(X_s)|^2\mathbf{1}_{\{s<\tau'\}}.
\end{align*}
Therefore, if we take $R$ with $\dis R< \frac{1}{2de^2a_2}$, we have
\begin{align}
d[\bi(X)^{\tau'},\bi(X)^{\tau'}]^c_t &\leq C_1dA_t,\label{Ageq} \\
d[\bi(X)^{\tau'},\bi(X)^{\tau'}]^d_t\mathbf{1}_{\{t<\tau' \}} &\leq C_1dB_t,\label{Bgeq}
\end{align}
where $\dis C_1=C_1(M,R,\alpha,\beta)=\left(-\frac{de a_2}{R} + \frac{e^{-1}}{2R^2}\right)^{-1}>0$. Thus, by \eqref{Ageq} and \eqref{Bgeq}, we have
\begin{align*}
[\bi(X),\bi(X)]^p_{\tau'-} &\leq \left( |\bi(X_0)|^2+C_1 (A_{\tau'}+B_{\tau'-})\right)^p\\
&\leq 2^{p-1}\beta^{2p} + 2^{p-1}C_1^p |Z_{\tau'-}-N_{\tau'-}|^p\\
&\leq 2^{p-1}\beta^{2p}+4^{p-1}C_1^p \left( (de)^p+\sup_{0\leq t < \infty}|N_{t\land \tau'}|^p \right).
\end{align*}
On the other hand, by Burkholder-Davis-Gundy's inequality, there exists $C_2=C_2(p)>0$ such that
\begin{align*}
\Ex \left[\left( \sup_{0\leq t <\infty}|N_{t\land \tau'}| \right)^p \right]\leq C_2\Ex \left[ [N,N]^{\frac{p}{2}}_{\tau'} \right].
\end{align*}
In addition, it holds that
\begin{align*}
[N,N]_t&=\sum_{i,j=1}^d \int_{0+}^t e^{Y^i_{s-}+Y^j_{s-}}D_k\bi^i(X_{s-}) D_l\bi^j(X_{s-})\, d[\bi^k(X),\bi^l(X)]_s\\
&\leq d^3a_1^2e^2[\bi(X),\bi(X)]_t.
\end{align*}
Therefore, we obtain
\begin{align*}
\Ex \left[ [\bi(X),\bi(X)]^p_{\tau'} \right]&\leq 2^{p-1}\left(\Ex \left[ [\bi(X),\bi(X)]^p_{\tau'-} \right] +\alpha^{2p} \right)\\
&\leq C_3+ C_4 \Ex \left[ [\bi(X),\bi(X)]_{\tau'}^{\frac{p}{2}} \right],
\end{align*}
where $C_3, C_4$ are constants which depend on $M,R,p,\alpha,\beta$. By setting
\begin{align*}
v&:=\Ex \left[ [\bi(X),\bi(X)]^p_{\tau'} \right]^{\frac{1}{2}},
\end{align*}
we have
\[
v^2\leq C_3+C_4v
\]
by Jensen's inequality. If we substitute $\sigma$ with
\begin{align}\label{sigman}
\sigma_n:=\inf \{ t\geq 0 \mid [\bi(X),\bi(X)]^p_t \geq n \},
\end{align}
then $v<\infty$ and we have
\[
v\leq \frac{C_4}{2}+\sqrt{C_3+\frac{C_4^2}{4}}.
\]
By letting $n$ tend to infinity, we obtain the desired estimate.
\end{proof}
The following lemma concerns Burkholder-Davis-Gundy's inequality for martingales on manifolds. In cases of continuous martingales on manifolds, a more sophisticated estimate was shown in Lemma 2.12 of \cite{AT98} by employing a suitable convex function on a small region. In our cases including discontinuous martingales, we need to employ another function defined on the external Euclidean space since the notion of martingales depends on the embedding.
\begin{lem}\label{Diffquad}
Let $\alpha,\beta>0$. Let $R>0$ be a constant constructed in \lref{Exquad}. Then for any $p\in [1,\infty)$, there exists $C=C(M,R,p,\alpha,\beta)$ such that
\[
\Ex \left[ [\bi(X)-\bi(Y),\bi(X)-\bi(Y)]^p_{t \land \tau} \right] \leq C \Ex \left[ \left( \sup_{0\leq s \leq t\land \tau}|\bi(X_s)-\bi(Y_s)| \right)^{2p} \right]^{\frac{1}{2}}
\]
for every $t\geq 0$, $(X, \zeta, \bp), (Y, \zeta', \bp) \in \mathcal{M}_{\alpha,\beta}(M,\bp)$ and stopping time $\tau$ satisfying \eqref{smallrange} and \eqref{smallrangeY}.
\end{lem}
\begin{proof}
Let $(X, \zeta, \bp), (Y, \zeta', \bp) \in \mathcal{M}_{\alpha,\beta}(M,\bp)$ and set $W=\bi(X)-\bi(X_0)$, $Z=\bi(Y)-\bi(Y_0)$. Let $h(w,z)=e^{|z-w|^2}-1$. Here we denote the canonical coordinate of $\mathbb{R}^d \times \mathbb{R}^d$ by $(w,z)=(w^1,\dots,w^d,z^1,\dots,z^d)$. For simplicity, we denote
\begin{align*}
\partial_i:=\frac{\partial }{\partial w^i},\ \partial_{\bar{i}}:=\frac{\partial }{\partial z^i}.
\end{align*}
Then by It\^o's formula, we have
\[
h(W_t,Z_t)-h(W_0,Z_0)=N_t+A^{(1)}_t+A^{(2)}_t+B^{(1)}_t+B^{(2)}_t,
\]
where
\begin{align*}
N_t&=\int_{0+}^t\partial_ih(W_{s-},Z_{s-}) \partial_k\bi^i(X_{s-})\, d\bi^k(X_s) + \int_{0+}^t\partial_{\bar{i}}h(W_{s-},Z_{s-})\partial_k\bi^i(Y_{s-})\, d\bi^k(Y_s),\\
A^{(1)}_t&=\frac{1}{2}\int_{0+}^t\partial_i\partial_jh(W_{s-},Z_{s-})\, d[W^i,W^j]^c_s +\int_{0+}^t\partial_i\partial_{\bar{j}}h(W_{s-},Z_{s-})\, d[W^i,Z^j]^c_s\\
&\h +\frac{1}{2}\int_{0+}^t \partial_{\bar{i}}\partial_{\bar{j}}h(W_{s-},Z_{s-})\, d[Z^i,Z^j]^c_s,\\
A^{(2)}_t&=\frac{1}{2}\int_{0+}^t \partial_ih(W_{s-},Z_{s-})\partial_k\partial_l \bi^i(X_{s-})\, d[\bi^k(X),\bi^l(X)]^c_s\\
&\h +\frac{1}{2}\int_{0+}^t \partial_{\bar{i}}h(W_{s-},Z_{s-})\partial_k\partial_l \bi^i(Y_{s-})\, d[\bi^k(Y),\bi^l(Y)]^c_s,\\
B^{(1)}_t&=\sum_{0<s\leq t}\left\{h(W_s,Z_s)-h(W_{s-},Z_{s-})-\partial_kh(W_{s-},Z_{s-})\Delta W^k_s-\partial_{\bar{k}}h(W_{s-},Z_{s-})\Delta Z^k_s \right\},\\
B^{(2)}_t&=\sum_{0<s\leq t}\partial_ih(W_{s-},Z_{s-})\left\{ \bi^i(X_s)-\bi^i(X_{s-})-\langle \nabla \bi^i(X_{s-}),\Delta \bi(X_s)\rangle  \right\}\\
&\h +\sum_{0<s\leq t}\partial_{\bar{i}}h(W_{s-},Z_{s-})\left\{ \bi^i(Y_s)-\bi^i(Y_{s-})-\langle \nabla \bi^i(Y_{s-}),\Delta \bi(Y_s)\rangle  \right\}.
\end{align*}
Note that
\begin{align*}
\partial_ih(w,z)&=-2(z^i-w^i)e^{|z-w|^2},\\
\partial_{\bar{i}}h(w,z)&=2(z^i-w^i)e^{|z-w|^2}
\end{align*}
and
\begin{align*}
\partial_i\partial_jh(w,z)&=2e^{|z-w|^2}\left\{ 2(z^i-w^i)(z^j-w^j)+\delta_{ij} \right\},\\
\partial_i\partial_{\bar{j}}h(w,z)&=-2e^{|z-w|^2}\left\{ 2(z^i-w^i)(z^j-w^j)+\delta_{ij} \right\},\\
\partial_{\bar{i}}\partial_{\bar{j}}h(w,z)&=2e^{|z-w|^2}\left\{ 2(z^i-w^i)(z^j-w^j)+\delta_{ij} \right\}.
\end{align*}
Therefore, the Hessian of $h$ can be written as
\begin{align*}
\mathrm{Hess}\, h=2e^{|z-w|^2}
\left[
\begin{array}{ll}
 H+I & -(H+I)\\
-(H+I) & H+I
\end{array}
\right],
\end{align*}
where $I=\mathrm{id}_{\mathbb{R}^d}$, $H=(H_{ij})_{1\leq i,j\leq d}$, 
\[
H_{ij}=2(z^i-w^i)(z^j-w^j).
\]
Thus for all $u,v\in \mathbb{R}^d$,
\begin{align*}
[^tu\ ^tv]\mathrm{Hess}\, h(w,z)
\left[
\begin{array}{cc}
u\\
v
\end{array}
\right] &= 2e^{|z-w|^2} \left(^t(u-v)H(u-v)+|u-v|^2\right)\\
&\geq |u-v|^2.
\end{align*}
We set $a_i:=a_i(R+\alpha+\beta)$, $i=1,2$. Then it holds that
\begin{align}
A^{(1)}_t &\geq \frac{1}{2} [W-Z,W-Z]^c_t,\label{A1geq}\\
\int_0^t \left| dA^{(2)}_s \right| &\leq a_2 \left([\bi(X),\bi(X)]_t^c+[\bi(Y),\bi(Y)]_t^c \right) \sup_{0\leq s< t}|\nabla h(W_s,Z_s)| \nonumber \\
&\leq a_2 b_1 \left([\bi(X),\bi(X)]_t^c+[\bi(Y),\bi(Y)]_t^c \right) \sup_{0\leq s< t}|W_s-Z_s|,\label{A2leq}\\
B^{(1)}_t &=\sum_{0<s\leq t} \int_0^1 [^t(\Delta W_s)\ ^t(\Delta Z_s)] \mathrm{Hess}(\theta)
\left[
\begin{array}{cc}
\Delta W_s\\
\Delta Z_s
\end{array}
\right] (1-\theta)\, d\theta \nonumber \\
&\geq \sum_{0<s\leq t}\frac{1}{2} |\Delta W_s-\Delta Z_s|^2, \label{B1geq}\\
\int_0^t\left| dB^{(2)}_s \right| &\leq a_2 \left([\bi(X),\bi(X)]^d_t+[\bi(Y),\bi(Y)]^d_t \right) \sup_{0\leq s< t}|\nabla h(W_s,Z_s)| \nonumber \\
&\leq a_2 b_1 \left([\bi(X),\bi(X)]^d_t+[\bi(Y),\bi(Y)]^d_t \right) \sup_{0\leq s< t}|W_s-Z_s|,\label{B2leq}
\end{align}
where 
\begin{align*}
\mathrm{Hess}(\theta)&=\mathrm{Hess}\, h ((W_{s-},Z_{s-})+\theta (\Delta W_s, \Delta Z_s)),\\
b_1&=2\sqrt{2}e^{4(R+\alpha)^2}.
\end{align*}
In addition, we have
\begin{align*}
[N^{\tau},N^{\tau}]_t&\leq 2\int_0^t |\partial_ih(W_{s-},Z_{s-})\partial_k \bi^i(X_{s-})|^2\, d[\bi^k(X),\bi^k(X)]_s\\
&\h +2\int_0^t |\partial_{\bar{i}}h(W_{s-},Z_{s-})\partial_k \bi^i(Y_{s-})|^2\, d[\bi^k(Y),\bi^k(Y)]_s\\
&\leq 2 a_1^2 \left( [\bi(X),\bi(X)]^{\tau}_t+[\bi(Y),\bi(Y)]^{\tau}_t \right) \sup_{0\leq s<t\land \tau}|\nabla h(W_s,Z_s)|^2\\
&\leq a_1^2 b_1^2 \left( [\bi(X),\bi(X)]^{\tau}_t+[\bi(Y),\bi(Y)]^{\tau}_t \right) \sup_{0\leq s<t\land \tau}|W_s-Z_s|^2.
\end{align*}
Thus by Burkholder-Davis-Gundy's inequality, there exists $C_1=C_1(p)$ such that
\begin{align}
\Ex \left[\left( \sup_{0\leq s \leq t} |N_{s\land \tau}|\right)^p \right] &\leq C_1 \Ex \left[ [N,N]^{\frac{p}{2}}_{t\land \tau} \right] \nonumber \\
&\leq C_1a_1^pb_1^p \Ex \left[  \left( [\bi(X),\bi(X)]^{\tau}_t+[\bi(Y),\bi(Y)]^{\tau}_t \right)^{\frac{p}{2}} \left( \sup_{0\leq s<t\land \tau}|W_s-Z_s|^2 \right)^{\frac{p}{2}} \right] \nonumber \\
&\leq 2^{\frac{p-1}{2}}C_1^{\frac{1}{2}} C_2^{\frac{1}{2}}a_1^pb_1^p \Ex \left[ \left( \sup_{0\leq s < t\land \tau}|W_s-Z_s|^2 \right)^p \right]^{\frac{1}{2}}, \label{Nleq}
\end{align}
where $C_2$ is the constant for $p$ in \lref{Exquad}. Therefore, by \eqref{A1geq}, \eqref{A2leq}, \eqref{B1geq}, \eqref{B2leq} and \eqref{Nleq}, it holds that
\begin{align*}
\Ex \left[ [W-Z,W-Z]^p_{t\land \tau}\right] &\leq 2^p \left[ \left( A^{(1)}_{t\land \tau}+B^{(1)}_{t\land \tau} \right)^p\right] \\
&= 2^p \Ex \left[ \left( h(W_{t\land \tau},Z_{t\land \tau})-h(W_0,Z_0)-N_{t\land \tau}-A^{(2)}_{t\land \tau}-B^{(2)}_{t\land \tau} \right)^p \right] \\
&\leq 2^{5p-4} \left( 2^p b_0^p \Ex \left[ \left( \sup_{0<s\leq t\land \tau}|W_s-Z_s| \right)^p \right]\right. \\
&\h +2^{\frac{p-1}{2}}C_1^{\frac{1}{2}}C_2^{\frac{1}{2}}a_1^pb_1^p \Ex \left[ \left( \sup_{0<s\leq t\land \tau}|W_s-Z_s|^2 \right)^p \right]^{\frac{1}{2}}\\
&\h +2^pa_2^pb_1^pC_1^{\frac{1}{2}}\Ex \left[ \left( \sup_{0<s\leq t\land \tau}|W_s-Z_s| \right)^{2p} \right]^{\frac{1}{2}}\\
&\h \left. +2^pa_2^pb_1^pC_1^{\frac{1}{2}}\Ex \left[ \left( \sup_{0<s\leq t\land \tau}|W_s-Z_s| \right)^{2p} \right]^{\frac{1}{2}} \right)\\
&\leq  C_3 \Ex \left[ \left( \sup_{0<s\leq t\land \tau}|W_s-Z_s| \right)^{2p} \right]^{\frac{1}{2}},
\end{align*}
where
\[
b_0=4(R+\alpha)e^{4(R+\alpha)^2}
\]
and $C_3$ is a constant which depends only on $M,p,R,\alpha,\beta$. Since it holds that
\begin{align*}
[\bi(X)-\bi(Y),\bi(X)-\bi(Y)]_{t\land \tau} &\leq 2[W-Z,W-Z]_{t\land \tau},\\
\sup_{0\leq s \leq t\land \tau} |W_s-Z_s| &\leq 2 \sup_{0\leq s \leq t\land \tau}|\bi(X_s)-\bi(Y_s)|,
\end{align*}
we obtain the desired estimate.
\end{proof}
\begin{proof}[Proof of \tref{smtopology}]
Let $(X,\zeta,\bp)$ and $(Y,\zeta',\bp)$ be $M$-valued martingales with the end point $\bp$ in $\mathcal{M}_{\alpha,\beta}(M,\bp)$. By It\^o's formula, we have
\begin{align*}
\bi^i(X_t)-\bi^i(Y_t)&=\bi^i(X_0)-\bi^i(Y_0)+N^{\iota^i}(X)_t-N^{\iota^i}(Y)_t+A^{\iota^i}(X)_t-A^{\iota^i}(Y)_t\\
&\h +B^{\iota^i}(X)_t-B^{\iota^i}(Y)_t,
\end{align*}
where
\begin{align*}
N^{\iota^i}(X)_t&:=\int_{0+}^t \langle \nabla \bi^i (X_{s-}),d\bi(X_s)\rangle,\\
A^{\iota^i}(X)_t&:=\frac{1}{2}\int_{0+}^t \nabla d\iota^i (X_{s-})\, d[X,X]^c_s,\\
B^{\iota^i}(X)_t&:=\sum_{0<s\leq t} \{ \bi^i(X_s) - \bi^i (X_{s-}) -\langle \nabla \bi^i(X_{s-}),\Delta \bi(X_s) \rangle \}
\end{align*}
and $N^{\iota^i}(Y), A^{\iota^i}(Y), B^{\iota^i}(Y)$ are defined in the same way. We set $a_i=a_i(R+\alpha+\beta),\ i=1,2,3$. Since
\begin{align*}
N^{\iota^i}(X)_t-N^{\iota^i}(Y)_t&=\int_{0+}^t\partial_j\bi^i(X_{s-})\, d(\bi^j(X_s)-\bi^j(Y_s))\\
&\h +\int_{0+}^t\left( \partial_j\bi^i(X_{s-})-\partial_j\bi^i(Y_{s-}) \right)\, d\bi^j(Y_s),
\end{align*}
we have
\begin{align}
\sum_{i=1}^d[N^{\iota^i}(X)-N^{\iota^i}(Y),N^{\iota^i}(X)-N^{\iota^i}(Y)]_t &\leq 2d a_1^2[\bi(X)-\bi(Y),\bi(X)-\bi(Y)]_t \nonumber\\[-2.5ex]
&\h +2\sum_{i=1}^d\sup_{0\leq s< t}|D \bi^i(X_s)-D \bi^i(Y_s)|^2 [\bi(Y),\bi(Y)]_t \nonumber\\
&\leq 2d a_1^2[\bi(X)-\bi(Y),\bi(X)-\bi(Y)]_t \nonumber\\
&\h +2d a_2^2[\bi(Y),\bi(Y)]_t \sup_{0\leq s< t}|\bi(X_s)-\bi(Y_s)|^2.\label{Niotaleq}
\end{align}
As for the continuous locally finite variation part, we have
\begin{align}
2\int_{0+}^t\left|d\left(A^{\iota^i}(X)-A^{\iota^i}(Y)\right)_s\right| &\leq \int_{0+}^t\left|\mathrm{Hess}\, \bi^i(X_{s-})\, d\left([\bi(X),\bi(X)]^c_s-[\bi(Y),\bi(Y)]^c_s\right)\right| \nonumber\\
&\h +\int_{0+}^t \left| \left( \mathrm{Hess}\, \bi^i(X_{s-})-\mathrm{Hess}\, \bi^i(Y_{s-}) \right)\, d[\bi(Y),\bi(Y)]^c_s\right| \nonumber\\
&\leq a_2\int_0^td\left|[\bi(X),\bi(X)]^c_s-[\bi(Y),\bi(Y)]^c_s\right| \nonumber\\
&\h + \sup_{0\leq s < t}\| \mathrm{Hess}\, \bi^i(X_s)-\mathrm{Hess}\, \bi^i(Y_s)\| [\bi(Y),\bi(Y)]^c_t \nonumber\\
&\leq a_2\left([\bi(X)+\bi(Y),\bi(X)+\bi(Y)]^c_t \cdot [\bi(X)-\bi(Y),\bi(X)-\bi(Y)]^c_t \right)^{\frac{1}{2}} \nonumber\\
&\h +a_3\sup_{0\leq s < t}|\bi(X_s)-\bi(Y_s)|[\bi(Y),\bi(Y)]^c_t. \label{Aiotaleq}
\end{align}
We can obtain a similar estimate for the discontinuous locally finite variation part. By Taylor's theorem, we have
\begin{align*}
B^{\iota^i}(X)_t-B^{\iota^i}(Y)_t&=\sum_{0<s\leq t} \int_0^1\left\{ ^t(\Delta \bi(X_s)) \mathrm{Hess}\, \bi^i(\bi(X_{s-})+\theta \Delta \bi (X_s)) \Delta \bi(X_s) \right.\\
&\h \left.- ^t(\Delta \bi(Y_s)) \mathrm{Hess}\, \bi^i(\bi(Y_{s-})+\theta \Delta \bi(Y_s)) \Delta \bi(Y_s) \right\} (1-\theta) \, d\theta.
\end{align*}
Since it holds that
\begin{align*}
&^t(\Delta \bi(X_s)) \mathrm{Hess}\, \bi^i(\bi(X_{s-})+\theta \Delta \bi (X_s))\Delta \bi(X_s)-^t(\Delta \bi(Y_s)) \mathrm{Hess}\, \bi^i(\bi(Y_{s-})+\theta \Delta \bi(Y_s))\Delta \bi(Y_s)\\
&=^t(\Delta \bi(X_s)) \left\{ \mathrm{Hess}\, \bi^i(\bi(X_{s-})+\theta \Delta \bi (X_s))-\mathrm{Hess}\, \bi^i(\bi(Y_{s-})+\theta \Delta \bi(Y_s))\right\}\Delta \bi(X_s)\\
&\h +^t(\Delta \bi(X_s)+\Delta \bi(Y_s))\mathrm{Hess}\, \bi^i(\bi(Y_{s-})+\theta \Delta \bi(Y_s))(\Delta \bi(X_s)-\Delta \bi(Y_s)),
\end{align*}
we have
\begin{align}
&\sum_{i=1}^d\int_{0+}^t\left|d\left(B^{\iota^i}(X)-B^{\iota^i}(Y)\right)_s\right| \nonumber\\
&\leq \sup_{0\leq s< t, \theta \in [0,1], i\in \{1,\dots,d \}}\| \mathrm{Hess}\, \bi^i(\bi(X_{s-})+\theta \Delta \bi(X_s))-\mathrm{Hess}\, \bi^i(\bi(Y_{s-})+\theta \Delta \bi(Y_s))\| \nonumber\\
&\h \cdot  \sum_{0<s\leq t} |\Delta \bi(X_s)|^2+a_2 \sum_{0<s\leq t}|\Delta(\bi(X)+\bi(Y))_s|\cdot |\Delta(\bi(X)-\bi(Y))_s| \nonumber\\
&\leq a_3\sup_{0\leq s< t, \theta \in [0,1]} |(\bi(X_{s-})+\theta \Delta \bi(X_s))-(\bi (Y_{s-})+\theta \Delta \bi (Y_s))| [\bi(X),\bi(X)]_t^d \nonumber\\
&\h +a_2\left([\bi(X)+\bi(Y),\bi(X)+\bi(Y)]_t^d \cdot [\bi(X)-\bi(Y),\bi(X)-\bi(Y)]_t^d\right)^{\frac{1}{2}}. \label{Biotaleq}
\end{align}
For the simplicity of notation, we set
\[
S(X,Y)_t=\sup_{0\leq s\leq t}|\bi (X_s)-\bi(Y_s)|.
\]
We note that
\begin{align*}
\sup_{0\leq s< t, \theta \in [0,1]} |(\bi(X_{s-})+\theta \Delta \bi(X_s))-(\bi (Y_{s-})+\theta \Delta \bi (Y_s))| \leq S(X,Y)_t\ \text{for all}\ t\geq 0\ \text{a.s.}
\end{align*}
Therefore, by Lemmas \ref{Exquad}, \ref{Diffquad} and \eqref{Niotaleq}, \eqref{Aiotaleq}, \eqref{Biotaleq}, for an arbitrary fixed $t>0$, we obtain
\begin{align*}\label{XYH1}
\sum_{i=1}^d\|\bi^i(X)^{t\land \tau}&-\bi^i(Y)^{t\land \tau} \|^p_{\mathbb{H}^p}\\[-2.5ex]
&\leq \sum_{i=1}^d\Ex \left[ |\bi^i(X_0)-\bi^i(Y_0)|^p\right. \\
&\h \left. + [N^{\iota^i}(X)^{t\land \tau}-N^{\iota^i}(Y)^{t\land \tau},N^{\iota^i}(X)^{t\land \tau}-N^{\iota^i}(Y)^{t\land \tau}]_{\infty}^{\frac{p}{2}}\right.\nonumber \\
&\h +\left( \int_0^{\infty} \left|d\left(A^{\iota^i}(X)^{t\land \tau}-A^{\iota^i}(Y)^{t\land \tau}\right)_s\right| \right)^p \nonumber \\
&\h \left.+\left( \int_0^{\infty} \left|d\left(B^{\iota^i}(X)^{t\land \tau}-B^{\iota^i}(Y)^{t\land \tau}\right)_s\right| \right)^p \right] \nonumber \\
&\leq \Ex \left[ S(X,Y)_{t\land \tau}^p \right] +2^p d^{\frac{p}{2}} a_1^p\Ex \left[ [\bi(X)-\bi(Y),\bi(X)-\bi(Y)]_{t\land \tau}^{\frac{p}{2}} \right] \\
&\h +2^pd^{\frac{p}{2}} a_2^p \Ex\left[ [\bi(Y),\bi(Y)]_{t\land \tau}^\frac{p}{2}S(X,Y)_{t\land \tau}^p \right]\\
&\h +2^{p-1}a_3^p \Ex \left[ \left( [\bi(Y),\bi(Y)]^c_{t\land \tau} S(X,Y)_{t\land \tau} \right)^p \right] \nonumber \\
&\h +2^{p-1}a_2^p\Ex \left[ \left([\bi(X)+\bi(Y),\bi(X)+\bi(Y)]^c_{t\land \tau}\cdot [\bi(X)-\bi(Y),\bi(X)-\bi(Y)]^c_{t\land \tau} \right)^{\frac{p}{2}} \right] \\
&\h +2^{p-1}a_3^p \Ex \left[ \left([\bi(X),\bi(X)]^d_{t\land \tau} S(X,Y)_{t\land \tau}\right)^p \right] \\
&\h +2^{p-1}a_2^p \Ex \left[ \left([\bi(X)+\bi(Y),\bi(X)+\bi(Y)]_{t\land \tau}^d \cdot [\bi(X)-\bi(Y),\bi(X)-\bi(Y)]_{t\land \tau}^d\right)^{\frac{p}{2}} \right] \\
&\leq \Ex \left[ S(X,Y)^p_{t\land \tau} \right]+2^pd^{\frac{p}{2}}a_1^pC_2(p)^{\frac{1}{2}} \Ex \left[ S(X,Y)_{t\land \tau}^{2p} \right]^{\frac{1}{4}} \\
&\h +2^pd^{\frac{p}{2}}a_2^pC_1(p)^{\frac{1}{2}}\Ex \left[ S(X,Y)_{t\land \tau}^{2p} \right]^{\frac{1}{2}} +2^{p}a_3^p C_1(2p)^{\frac{1}{2}} \Ex \left[ S(X,Y)_{t\land \tau}^{2p} \right]^{\frac{1}{2}}\\
&\h +2^{p+1} a_2^p C_1(p)^{\frac{1}{2}} C_2(p)^{\frac{1}{2}}\Ex \left[ S(X,Y)_{t\land \tau}^{2p} \right]^{\frac{1}{4}},
\end{align*}
where $C_i(p)$, $i=1,2$ are the constants in \lsref{Exquad}, \ref{Diffquad} for the parameter $p$, respectively. Since $S(X,Y)\leq 2(R+\alpha+\beta)$, there exists a constant $C=C(M,R,p,\alpha,\beta)$ such that
\[
\left( \sum_{i=1}^d\|\bi^i(X)^{t\land \tau}-\bi^i(Y)^{t\land \tau} \|^p_{\mathbb{H}^p}\right)^{\frac{1}{p}} \leq C \Ex \left[ S(X,Y)_{t\land \tau}^{2p} \right]^{\frac{1}{4p}}.
\]
This completes the proof.
\end{proof}
Next we will show \cref{equitopo}. In the statement of \cref{equitopo}, the existence of $\zeta$ is elementary.
\begin{lem}\label{killingX}
Let $\{(X^n, \zeta^n,\bp)\}_{n\in \mathbb{N}}$ be a sequence of $M$-valued semimartingales with an end point $\bp$. Assume that $X^n$ converges to a \cl process $X$ uniformly in every compact interval almost surely on $\mathbb{R}^d$. Then there exists a stopping time $\zeta$ satisfying \eqref{liminfzeta} such that for almost all $\omega \in \Omega$, $X_t(\omega)\in M$ for $t<\zeta(\omega)$ and $X_t(\omega) =\bp$ for $t\geq \zeta(\omega)$. 
\end{lem}
\begin{proof}
Since $M$ is closed in $\mathbb{R}^d$ and $\bp \notin M$, we can define
\[
\zeta:=\inf \{ t\geq 0 \mid X_t \notin M\}.
\]
In addition, since the limit process $X$ is adapted by assumption, the random time $\zeta$ is a stopping time.
Next, we fix $T>0$. For almost all $\omega \in \{ \zeta \leq T\}$, there exists $N_1(\omega)\in \mathbb{N}$ such that $X^n_{\zeta(\omega)}(\omega)=\bp$ for all $n \geq N_1(\omega)$. This means $\zeta^n(\omega)\leq \zeta(\omega)$. On the other hand, there exists $N_2(\omega)\in \mathbb{N}$ such that $X^n_{\zeta(\omega)-}(\omega)\in M$ for all $n\geq N_2(\omega)$. Thus $\zeta^n(\omega)\geq \zeta(\omega)$. Therefore, we have $\zeta^n(\omega)=\zeta(\omega)$ for $n \geq N_1(\omega) \lor N_2(\omega)$. Finally, for almost all $\omega \in \{ \zeta > T \}$, there exists $N_3(\omega) \in \mathbb{N}$ such that for all $n\geq N_3(\omega)$ and $t \in [0,T]$, $X^n_t(\omega) \in M$. This means $\zeta^n(\omega)\land T=\zeta(\omega)\land T=T$ for $n \geq N_3(\omega)$. Therefore we have \eqref{liminfzeta}.
\end{proof}
In the proof of \cref{equitopo}, we will restrict the scope to sequences of martingales valued in a small open subset. The technique by the following lemma might be well known and frequently used in the study of continuous martingales on manifolds (see the proof of Theorem 4.43 of \cite{Emery89}) but we confirm it here in order to apply it to martingales with jumps.
\begin{lem}\label{localize}
Let $\{(X^n,\zeta^n,\bp)\}_{n\in \mathbb{N}}$ be a sequence of $M$-valued semimartingales with an end point $\bp$.
We further assume that $\{X^n\}_{n\in \mathbb{N}}$ converges to an $\mathbb{R}^d$-valued \cl process $X$ uniformly on every compact interval almost surely. Let $\zeta$ be a killing time of $X$ defined in the proof of \lref{killingX}. Let $\{U_j\}_{j\in \mathbb{N}}$ be a countable open covering of $M$ which consists of geodesic balls. Let $\{ j(i)\}_{i\in \mathbb{N}}$ be a sequence of positive integers in which every positive integer appears infinitely many times. Then for any countable geodesic balls $\{V_j\}_{j\in \mathbb{N}}$ covering $M$ with $\overline{U_j} \Subset V_j$ for each $j$, there exist non-decreasing sequences of stopping times $\{\tau_i\}_{i=0}^{\infty}$, $\{ \sigma^{(i)}_{m}\}_{m=1}^{\infty}$ and $\{ \sigma_m\}_{m=1}^{\infty}$ satisfying the following:
\begin{itemize}
\item[(i)]$\tau_0=0$ and
\[
\Prob \left( \bigcap_{k=1}^{\infty} \liminf_{i\to \infty}\{ \tau_i \land k=\zeta \land k \}\right)=1.
\]
\item[(ii)]For each $i,m$, $\tau_i \leq \sigma_m^{(i)} \leq \tau_{i+1}$ and it holds that
\[
\Prob \left( \bigcap_{k=1}^{\infty} \liminf_{m\to \infty}\{ \sigma_m^{(i)} \land k=\tau_{i+1} \land k \}\right)=1.
\]
\item[(iii)]For each $m$, $\zeta \leq \sigma_m$ and it holds that
\begin{align*}
\Prob \left(\liminf_{m\to \infty}\{\sigma_m=\infty \}\right)=1.
\end{align*}
\item[(iv)]For each $m,n\in \mathbb{N}$ with $n\geq m$, it holds that
\begin{align*}
X^n_t(\omega) \in V_{j(i)}\ &\text{for a.e.}\ \omega \in \{\tau_i < \sigma_m^{(i)} \}\ \text{and}\ t\in [\tau_i(\omega), \sigma_m^{(i)}(\omega)),\\
X^n_t(\omega) =\bp \ &\text{for a.e.}\ \omega \in \{\zeta < \sigma_m \}\ \text{and}\ t\in [\zeta(\omega), \sigma_m(\omega)).
\end{align*}
\end{itemize}
\end{lem}
\begin{proof}
Note that $X_t \in M$ for $t<\zeta$ and $X_t=\bp$ for $t\geq \zeta$ in our setting. We define a sequence of stopping times $\{ \tau_i\}_{i\in \mathbb{N}}$ inductively by
\begin{align*}
\tau_0&=0,\\
\tau_{i+1}&=\inf \{t\geq \tau_i \mid X_t \notin U_{j(i)} \}.
\end{align*}
Then for each $i$, $X_t$ stays in $\overline{U}_{j(i)}$ on $[\tau_i,\tau_{i+1})$ if $\tau_i <\tau_{i+1}$.
We note that $\dis \lim_{i\to \infty}\tau_i=\zeta$ a.s. since every positive integer appears infinitely many times in the sequence $\{j(i)\}_{i\in \mathbb{N}}$. In addition, on the event $\{\zeta < \infty\}$, there exists $i^* \in \mathbb{N}$ such that $X_{\zeta-}\in U_{j(i^*)}$ and we have $\tau_{i^*}=\zeta$. Therefore, the sequence $\{\tau_i\}_{i\in \mathbb{N}}$ satisfies (i). We further set
\begin{align*}
\tilde \sigma_m^{(i)}&:=\inf \{t\geq \tau_i \mid X^m\notin V_{j(i)} \},\\
\sigma_m^{(i)}&:=\left( \inf_{n\geq m}\tilde \sigma_n^{(i)} \right) \land \tau_{i+1}.
\end{align*}
Then for a.s. $\omega \in \Omega$ with $\tau_i(\omega)<\tau_{i+1}(\omega)<\infty$, there exists $N(\omega)\in \mathbb{N}$ such that $m\geq N(\omega)$ yields
\[
\sup_{\tau_i(\omega) \leq s < \tau_{i+1}(\omega)} |\bi(X_s(\omega))-\bi(X_s^m(\omega))| \leq \delta_i,
\]
where $\dis \delta_i=\inf_{(x,y)\in \partial U_{j(i)}\times \partial V_{j(i)}}|\bi(y)-\bi(x)|$. This means $\sigma_m^{(i)}(\omega)= \tau_{i+1}(\omega)$ if $m\geq N(\omega)$ since $X_s(\omega)\in U_{j(i)}$ for $s\in [\tau_i(\omega),\tau_{i+1}(\omega))$. On the other hand, for a.s. $\omega \in \{\tau_{i+1}=\infty\}$ and each $k\in \mathbb{N}$, there exists $N(\omega)\in \mathbb{N}$ such that $m\geq N(\omega)$ yields
\[
\sup_{\tau_i(\omega) \leq s \leq k} |\bi(X_s(\omega))-\bi(X_s^m(\omega))| \leq \delta_i.
\]
Thus $\dis \lim_{m\to \infty} \sigma_m^{(i)}(\omega)=\infty$. Therefore $\{\sigma_m^{(i)}\}_{m\in \mathbb{N}}$ satisfies (ii). Next we fix $m\in \mathbb{N}$ and $n\geq m$. Then for $\omega \in \Omega$ with $\tau_i(\omega) < \sigma_m^{(i)}(\omega)$, it holds that $\tau_i(\omega)<\tilde \sigma_n^{(i)}(\omega)$. This yields
\begin{align*}
X^n_t(\omega)\in V_{j(i)}
\end{align*}
for $t\in [\tau_i(\omega),\sigma_m^{(i)}(\omega))$.
Next we set
\begin{align*}
\tilde \sigma_m&:=\inf \{ t\geq \zeta \mid X^m_t \neq p \},\\
\sigma_m&:=\inf_{n\geq m}\tilde \sigma_n.
\end{align*}
Then
\begin{align*}
\tilde \sigma_m(\omega)=
\begin{cases}
\zeta(\omega)\ &\text{on}\ \{\zeta < \zeta^m, \zeta < \infty\}\\
\infty \ &\text{on}\ \{\zeta^m \leq \zeta\}\cup \{\zeta=\infty\}.
\end{cases}
\end{align*}
By assumption, for a.s. $\omega \in \Omega$ with $\zeta(\omega)<\infty$, $\dis \lim_{n\to \infty}X^n_{\zeta(\omega)}(\omega)=X_{\zeta(\omega)}(\omega)=\bp$. This means that $X^n_{\zeta(\omega)}(\omega)=\bp$ for sufficiently large $n$ since $M$ is closed in $\mathbb{R}^d$ and $\bp \notin M$. Therefore, (iii) holds. Finally, we fix $m\in \mathbb{N}$ and take $\omega \in \{\zeta<\sigma_m\}$, $t\in [\zeta(\omega),\sigma_m(\omega))$ and $n \geq m$. Then $X^n_t(\omega)=\bp$ since $t\in [\zeta(\omega),\tilde \sigma_n(\omega))$. Thus we have (iv).
\end{proof}

\begin{lem}\label{smlocalize}
Let $X$ and $\{X^n\}_{n\in \mathbb{N}}$ be semimartingales. Let $\{\tau_i\}_{i\in \mathbb{N}}$ be a non-decreasing sequence of stopping times with $\tau_0=0$ and $\dis \lim_{i\to \infty}\tau_i = \infty$ a.s. Set for each $i\in \mathbb{N}$,
\begin{align*}
X^{n,(i)}_t&:=X^n_{(\tau_i+t)\land \tau_{i+1}},\\
X^{(i)}_t&:=X_{(\tau_i+t)\land \tau_{i+1}}.
\end{align*}
We assume the following:
\begin{itemize}
\item[(i)]For each $i$, there exists a non-decreasing sequence of $\{\mathcal{F}_{\tau_i+t}\}_{t\geq 0}$-stopping times $\sigma^{(i)}_m$ such that $\sigma^{(i)}_m \geq \tau_i$ and for each $T>0$,
\[
\Prob\left(\liminf_{m \to \infty}\{ \sigma^{(i)}_m \land T= \tau_{i+1} \land T \} \right)=1.
\]
\item[(ii)]For each fixed $m \in \mathbb{N}$, the sequence of $\{\cF_{\tau_i+t}\}_{t\geq 0}$-semimartingale $\{X^{n,(i),\theta^{(i)}_m}\}_{n \in \mathbb{N}}$ converges to $X^{(i),\theta_m^{(i)}}$ in the semimartingale topology as $n\to \infty$, where $\theta_m^{(i)}=\sigma_m^{(i)}-\tau_i$.
\end{itemize}
Then the sequence $\{X^n\}_{n\in \mathbb{N}}$ converges to $X$ in the semimartingale topology.
\end{lem}
\begin{proof}
We set $Y^n:=X^n-X$ and $Y^{n,(i)}=X^{n,(i)}-X^{(i)}$. It suffices to show that $\dis \lim_{n\to \infty}\hat r(Y^{n,\tau_i})=0$ for each $i$. Note that any $\{\cF_{\tau_i + t}\}_{t\geq 0}$-predictable process can be written in the form of $H_{\tau_i +\cdot}$ for some $H \in \cP$.
Therefore, it holds that
\begin{align*}
\hat r(Y^{n,(i)})&=\sup_{H\in \mathcal{P}_1}\sum_{k=1}^{\infty} \frac{1}{2^k} \Ex \left[ 1\land \sup_{0\leq t \leq k}\left| \int_0^t H_{(\tau_i+\cdot)}\, dY^{n,(i)}\right| \right]\\
&=\sup_{H\in \mathcal{P}_1}\sum_{k=1}^{\infty}\frac{1}{2^k} \left( \Ex \left[ 1\land \sup_{0\leq t \leq k}\left| \int_0^t H_{(\tau_i+\cdot)}\, dY^{n,(i),\theta_m^{(i)}}\right| ; \theta_m^{(i)}\geq k\land (\tau_{i+1}-\tau_i) \right] \right. \\
&\h \left.+  \Ex \left[ 1\land \sup_{0\leq t \leq k}\left| \int_0^t H_{(\tau_i+\cdot)}\, dY^{n,(i)}\right| ; \theta_m^{(i)} <k\land (\tau_{i+1}-\tau_i) \right]\right) \\
&\leq \hat r(Y^{n,i,\theta_m^{(i)}})+\sum_{k=1}^{\infty} \frac{1}{2^k}\Prob (\theta_m^{(i)}<k\land (\tau_{i+1}-\tau_i)),
\end{align*}
where we used
\[
Y^{n,(i)}_t=Y^{n,(i),\theta_m^{(i)}}_t
\]
for all $t\leq k$ on the event $\{\theta_m^{(i)} \geq k \land (\tau_{i+1}-\tau_i) \}$ in the second equality.
Therefore,
\[
\limsup_{n\to \infty} \hat r(Y^{n,(i)})\leq \sum_{k=1}^{\infty} \frac{1}{2^k}\Prob (\theta_m^{(i)}<k\land (\tau_{i+1}-\tau_i))
\]
for all $m$ by assumption (ii). By letting $m$ tend to $\infty$ and applying the bounded convergence theorem, we obtain $\dis \lim_{n\to \infty} \hat r(Y^{n,(i)})=0$ for each $i$ by assumption (i). Since for each $i$, $k$ and $H\in \mathcal{P}_1$,
\begin{align*}
\Ex \left[1\land \sup_{0\leq s \leq k} \left| \int_0^s H\, dY^{n,\tau_i}\right| \right] \leq \sum_{j=0}^{i-1}\Ex \left[1\land \sup_{0\leq s \leq k} \left| \int_0^s H_{\tau_j+s}\, dY_s^{n,(j)}\right|\ ; \tau_j<\infty \right],
\end{align*}
it holds that
\begin{align*}
\hat r(Y^{n,\tau_i})\leq \sum_{j=0}^{i-1}\hat r(Y^{n,(j)}) \to 0
\end{align*}
as $n\to \infty$ for each $i$. This yields $\hat r(Y^n)\to 0$ as $n \to \infty$.
\end{proof}

\begin{proof}[Proof of \cref{equitopo}]
Let $\{(X^n,\zeta^n,\bp)\}_{n\in \mathbb{N}}$ be a sequence in $\mathcal{M}_{\alpha,\beta}(M, \bp)$ converging to a process $X$ in the u.c.p. topology. Assume that an arbitrary subsequence is given. By extracting a further subsequence, we can assume that $\{X^n\}_{n\in \mathbb{N}}$ converges to $X$ on each compact interval almost surely, where we continue to use the index $n$ for the subsequence to simplify the notation. Let $\zeta$ be a stopping time constructed in \lref{killingX}. For $x\in M$, let $R_x=R(M,\alpha,|x|+1)$ be the radius obtained in \tref{smtopology} for $\beta=|x|+1$. Take geodesic balls $U_x, V_x$ with $U_x \Subset V_x$ so small that they satisfy
\[
\sup_{y\in V_x}|y| \leq |x|+1,\ \sup_{y_1,y_2 \in V_x}|y_1-y_2|< R_x.
\]
Since $M$ is second-countable, we can take countable points $\{x_j\}_{j\in \mathbb{N}}$ such that the open sets $U_j:=U_{x_j}$ cover $M$. Then we take a sequence $\{j(i)\}_{i\in \mathbb{N}}$ and construct sequences of stopping times $\{\tau_i\}_{i=0}^{\infty}$, $\{\sigma_m^{(i)}\}_{m=1}^{\infty}$ and $\{ \sigma_m\}_{m=1}^{\infty}$ as in \lref{localize} and set $\theta_m^{(i)}=\sigma_m^{(i)}-\tau_i$. We set $\tau_{\infty}=\zeta$, $\sigma_m^{(\infty)}=\sigma_m$ for the convenience of the notation. We set $X^{n,(i),\theta_m^{(i)}}$ as in \lref{smlocalize} for $m\in \mathbb{N}$, $i\in \mathbb{N} \cup \{ \infty \}$. We further set
\begin{align*}
\zeta_m^{n,(i)}:=
\begin{cases}
(\zeta_m^n-\tau_i)\lor 0\ &\text{on}\ \{ \zeta^n \leq \sigma_m^{(i)} \},\\
\infty \ &\text{on}\ \{ \zeta^n > \sigma_m^{(i)} \}.
\end{cases}
\end{align*}
For each $i,m$, we set $\Lambda_{i,m}:=\{ \tau_i < \sigma_m^{(i)} \}$. We assume $\Prob(\Lambda_{i,m})\in (0,1)$. Then for each $n,i,m$ with $n\geq m$,
\[
|X^{n,(i),\theta_m^{i}}_0|=|X^n_{\tau_i}| \leq |x_{j(i)}|+1\ \text{on}\ \Lambda_{i,m}
\]
by the construction of $V_j:=V_{x_j}$. Therefore, the triplet $(X^{n,(i),\theta_m^{(i)}}, \zeta_m^{n,(i)} ,\bp)$ belongs to $\mathcal{M}_{\alpha,|x_{j(i)}|+1}(M,\bp)$ with respect to $\Prob(\cdot \mid \Lambda_{i,m})$. Moreover, for fixed $i,m$ and each $n\geq m$, the process $X^{n,(i),\theta_m^{(i)}}$ stays in $V_{j(i)}$ on $[0,\theta_m^{(i)})$ $\Prob(\cdot \mid \Lambda_{i,m})$-a.s. Therefore, for any fixed $i,m$ and $t>0$, we have
\begin{align*}
\lim_{n,n'\to \infty}\sum_{k=1}^d\left\|\bi^k\left(X^{n,(i),t\land \theta_m^{(i)}} \right)-\bi^k \left(X^{n',(i),t\land \theta_m^{(i)}}\right) \right\|_{\mathbb{H}^2(\Lambda_{i,m})}=0
\end{align*}
by \tref{smtopology}, where $\mathbb{H}^2(\Lambda_{i,m})$ is the $\mathbb{H}^2$-space under $\Prob(\cdot \mid \Lambda_{i,m})$. On the other hand, $X^{n,(i),t\land \theta_m^{(i)}}$ is $\| \cdot \|_{\mathbb{H}^2(\Lambda_{i,j}^c)}$-Cauchy since $X^{n,(i),t\land \theta_m^{(i)}}$ does not move from $X^n_{\tau_i}$ on $\Lambda_{i,m}^c$. Therefore, $\{X^{n,(i),t\land \theta_m^{(i)}}\}_{n\in \mathbb{N}}$ is a Cauchy sequence in $\| \cdot \|_{\mathbb{H}^2}$. In the case of $\Prob(\Lambda_{i,m})=0$ or $1$, we can deduce the same consequence. Therefore, there exists a further subsequence which converges to $X$ in the semimartingale topology by \rref{Hpconv} and the argument of \lref{smlocalize} with the terminal interval $[\zeta,\infty)$ indexed by $i=\infty$. Therefore, the subsequence converges to $X$ in the semimartingale topology under $\Prob$. This holds for any subsequence of the original sequence of semimartingales $\{X^n\}_{n\in \mathbb{N}}$. Therefore, the original sequence also converges to $X$ in the semimartingale topology. Finally, since for any bounded vector field $\cX$ on $M$, the stochastic integral $\dis \int \langle \cX(X^n_-), dX^n\rangle$ is a local martingale with jumps bounded uniformly in $n$. Thus by the convergence in the semimartingale topology, we have
\[
\int \langle \cX(X^n_-), dX^n\rangle \to \int \langle \cX(X_-), dX \rangle
\]
in u.c.p. and the limit is a local martingale. Therefore, $(X,\zeta,\bp)$ is a martingale on $M$ with the end point $\bp$.
\end{proof}

\section{Proof of \tref{harmonicseq}}\label{PfofSeqharm}
To begin with, we recall some notions regarding Dirichlet forms following \cite{ChenFuku} and \cite{FOT}.

Let $(\mathcal{E},\mathcal{F})$ be a Dirichlet form on $L^2(E;\bdm)$, i.e. $\cF$ is a dense subspace of $L^2(E;\bdm)$, $\cE \colon \cF \times \cF \to \mathbb{R}$ is a symmetric bilinear Markovian form and $(\cF,\cE_1)$ is a Hilbert space. Here the inner product $\cE_1$ on $\cF$ is defined by
\[
\cE_1(u,v)=\cE(u,v) + \langle u, v \rangle_{L^2(E;\bdm)}
\]
A Dirichlet form $(\cE,\cF)$ is said to be regular if $C_0(E)\cap \cF$ is $\cE_1$-dense in $\cF$ and $\|\cdot \|_{\infty}$-dense in $C_0(E)$, where $C_0(E)$ is a space of continuous functions with compact supports. For an open set $D\subset E$, we set
\begin{align*}
\begin{cases}
\cF^D&:=\{ u \in \cF \mid \tilde u=0\ \text{q.e.}\ \text{on}\ E \backslash D \},\\
\cE^D(u,v)&:=\cE(u,v),\ u,v \in \cF^D,
\end{cases}
\end{align*}
where $\tilde{u}$ is a quasi-continuous modification of $u$ and ``q.e." stands for ``except for a zero capacity set with respect to $\cE_1$".
Then $(\cE^D, \cF^D)$ is also a regular Dirichlet form on $L^2(D;\bdm)$. The $M$-valued local Dirichlet space $\cF^D_{loc}(M)$ on $D$ is defined as follows: We say an $\bdm$-measurable function $u\colon E \to \mathbb{R}$ is in the space $\cF^D_{loc}$ if for any relatively compact open subset $D_1 \Subset D$, there exists $u' \in \cF^D$ such that
\[
u=u'\ \bdm \text{-a.e. on}\ D_1.
\]
Then, we set
\begin{align*}
\cF^D_{loc}(\mathbb{R}^d)&:=\{ u=(u^1,\dots,u^d) \colon E \to \bbR^d \mid u^i \in \cF^D_{loc}\ \text{for each}\ i=1,\dots,d \},\\
\cF^D_{loc}(M)&:=\{ u\in \cF^D_{loc}(\mathbb{R}^d)\mid u(z) \in M\ \text{a.e.}\ z\in E \}.
\end{align*}
It is known that for a regular Dirichlet form $(\cE,\cF)$ there exists an $\bdm$-symmetric Hunt process $\bZ=(\Omega, \{Z_t\}_{t\geq 0}, \zeta, \{\Prob_z\}_{z\in E})$ such that its transition function $\{p_t\}_{t\geq 0}$ satisfies
\[
\lim_{t\to 0}\frac{1}{t}\langle u-p_tu,u \rangle_{L^2(E;\bdm)} \nearrow \cE(u,u)
\]
for each $u \in \mathcal{F}$ bounded and Borel measurable, where $\Omega$ is a measurable space with some $\sigma$-algebra, $\zeta$ is the killing time meaning $Z$ is trapped at $\partial$ after $\zeta$, and $\Prob_z$ is a probability measure on $\Omega$ meaning the law of the process $Z$ starting at $z\in E$. Let $\{\cF^0_t\}_{t=0}^{\infty}$ be the filtration generated by $Z$, i.e.
\begin{align*}
\mathcal{F}^0_{\infty}&=\sigma(Z_s;\ s< \infty),\\
\mathcal{F}^0_t&=\sigma(Z_s;\ s\leq t).
\end{align*}
For a $\sigma$-finite measure $\mu$, we denote the $\mathbb{P}_{\mu}$-completion of $\mathcal{F}^0_{\infty}$ by $ \mathcal{F}^{\mu}_{\infty}$, where $\Prob_{\mu}$ is a $\sigma$-finite measure defined by
\[
\mathbb{P}_{\mu}(\Lambda)=\int_{E_{\partial}}\mathbb{P}_z(\Lambda)\, \mu(dz),\ \Lambda \in \mathcal{F}^0_{\infty}.
\]
Let $\mathcal{N}_{\mu}$ be the family of all $\mathbb{P}_{\mu}$-null sets in $\mathcal{F}^{\mu}_{\infty}$.
We set
\[
\mathcal{F}^{\mu}_t=\sigma(\mathcal{F}^0_t,\mathcal{N}_{\mu}).
\]
Denote the set of probability measures on $E_{\partial}$ by $\mathcal{P}(E_{\partial})$ and let
\[
\mathcal{F}^Z_t=\bigcap_{\mu \in \mathcal{P}(E_{\partial})} \mathcal{F}^{\mu}_t,\ t\in [0,\infty].
\]

We denote the set of positive smooth measures with respect to $(\cE,\cF)$ by $\sS$. Then for each positive continuous additive functional (PCAF) $A$, there exists $\mu \in \sS$ such that
\[
\lim_{t \to 0}\frac{1}{t}\Ex_{\bdm}\left[ \int_0^t v(Z_s)\, dA_s \right]=\int_E v\, d\mu \ \text{for each}\ v \in \mathcal{B}_+(E)
\]
and this is a one-to-one correspondence between PCAF's and positive smooth measures.
A subclass $\sS_0 \subset \sS$ is defined by the set of positive Borel measures $\mu$ such that $\mu$ is a positive Radon measure on $E$ and the functional $v \mapsto \int_E v \, d\mu$ defined on $C_0(E) \cap \cF$ is bounded with respect to the $\cE_1$-norm. Then for each $\mu \in \sS_0$, there exists $U_1\mu \in \cF$ such that
\[
\cE_1(U_1\mu, v)=\int_E \tilde{v}\, d\mu \ \text{for all}\ v\in \cF,
\]
which is called the $1$-potential of $\mu$.
Furthermore, a subclass $\sS_{00} \subset \sS$ is defined by
\[
\sS_{00}:=\{ \mu \in \sS_0 \mid \mu (E)<\infty,\ \| U_1\mu\|_{L^{\infty}(E)}<\infty \}.
\]
For an arbitrary subset $B \subset E$, we set
\[
\cL_B:=\{ v \in \cF \mid v\geq 1\ \text{q.e. on}\ B\}.
\]
Then if $\cL_B \neq \emptyset$, by Theorem 2.3.1 of \cite{ChenFuku}, there exists a unique minimizer $\min_{v \in \cL_B}\cE_1(v,v)$ which we denote by $e_B$ and call the $1$-equilibrium potential of $B$. Moreover, $e_B$ is a $1$-potential of some measure $\nu_B \in \sS_0$.

For $v\in \cF$, there exists a martingale additive functional (MAF) of finite energy $M^{[v]}$ and a continuous additive functional (CAF) of zero energy $N^{[v]}$ such that
\[
\tilde v(Z_t)-\tilde v(Z_0)=M^{[v]}_t + N^{[v]}_t,\ t\geq 0,\ \Prob_z \text{-a.s. for q.e.}\ z\in E
\]
by Theorem 5.2.2 of \cite{FOT}. We denote the Revuz measure associated with the PCAF $\langle M^{[v_1]},M^{[v_2]} \rangle$ by $\mu_{\langle v_1, v_2 \rangle}$ for $v_1,v_2 \in \cF$. In particular, we denote $\mu_{\langle v, v \rangle}=\mu_{\langle v\rangle}$ for $v\in \cF$. For a vector valued function $v=(v^1,\dots,v^d)$ with $v^i\in \cF$ for $i=1,\dots,d$, we set
\[
\mu_{\langle v \rangle}=\sum_{i=1}^d\mu_{\langle v^i \rangle}.
\]
Let $(N,H)$ be a L\'evy system of $\bZ$, that is, $N$ is a kernel from $E$ to $E_{\partial}$, $H$ is a PCAF of $Z$ satisfying
\begin{align}\label{LevySys}
\Ex_z \left[ \sum_{0<s\leq t}Y_sF(Z_{s-},Z_s) \right]=\Ex_z \left[ \int_0^tY_s \int_{E_{\partial}}F(Z_{s},w)\, N(Z_s,dw)dH_s \right]
\end{align}
for every nonnegative predictable process $Y$ and nonnegative Borel function $F$ on $E_{\partial}\times E_{\partial}$ satisfying $F(w,w)=0$. Let $\mu_H$ be the Revuz measure of $H$. Then for the jump measure $J$ and killing measure $k$ defined by
\begin{align*}
J(dzdw)&=N(z,dw)\mu_H(dz),\\
k(dz)&=N(z,\{\partial \})\mu_H(dz),
\end{align*}
we have the following Beurling-Deny decomposition
\begin{align}
\cE(v_1,v_2)&=\frac{1}{2}\mu_{\langle v_1,v_2 \rangle}^c(E)+\frac{1}{2}\int_{E \times E \backslash \mathrm{diag}(E)}(\tilde{v_1}(z)-\tilde{v_1}(w))(\tilde{v_2}(z)-\tilde{v_2}(w))\, J(dzdw)\nonumber \\
&\h +\int_E \tilde{v_1}\tilde{v_2} \, dk \label{BDdecom}
\end{align}
for $v_1,v_2 \in \cF$, where $\mu_{\langle v_1,v_2 \rangle}^c$ is the signed Revuz measure of $\langle M^{[v_1],c},M^{[v_2],c} \rangle$.

\begin{lem}\label{estE}
Let $\phi \in C_0(D)\cap \cF$ and set $K=\mathrm{supp}[\phi]$. Then there exists a constant $C>0$ which depends only on $\phi$ such that for all $v \in \cF^D \cap L^{\infty}(E)$, 
\begin{align*}
\cE(\phi v, \phi v) \leq C \left( \mu_{\langle v \rangle}^c(K)+\int_{K \times K}|v(z)-v(w)|^2\, J(dzdw) + \| v \|_{L^{\infty}(K)}^2 \right).
\end{align*}
\end{lem}
\begin{proof}
By the Beurling-Deny decomposition,
\begin{align*}
\cE(\phi v, \phi v)=\frac{1}{2}\mu_{\langle \phi v \rangle}^c(D)+\frac{1}{2}\int_{E \times E}|\phi v(z)-\phi v(w)|^2\, J(dzdw)+\int_D |\phi v|^2\, dk.
\end{align*}
For the continuous part, since
\begin{align*}
d\mu_{\langle \phi v \rangle}^c=\phi^2 d\mu_{\langle v \rangle}^c + 2\phi v d\mu_{\langle \phi,  v \rangle}^c + v^2 d\mu_{\langle \phi \rangle}^c
\end{align*}
by Lemma 4.3.6 of \cite{ChenFuku}, we have
\[
\mu_{\langle \phi v \rangle}^c(D)=\mu_{\langle \phi v \rangle}^c (K) \leq 2\left(\| \phi \|_{L^{\infty}(D)}^2 \mu_{\langle v \rangle}^c(K) + \| v \|_{L^{\infty}(K)}^2 \mu_{\langle \phi \rangle}^c(K) \right).
\]
For the jump part, we have
\begin{align*}
\frac{1}{2}\int_{E \times E}|\phi v(z)-\phi v(w)|^2\, J(dzdw)&= \frac{1}{2}\int_{K \times K}|\phi v(z)-\phi v(w)|^2\, J(dzdw)\\
&\h + \int_{K \times K^c}|\phi v(z)|^2\, J(dzdw)\\
&\leq \|\phi\|_{L^{\infty}(K)}^2\int_{K \times K}|v(z)-v(w)|^2\, J(dzdw)\\
&\h + \|v \|_{L^{\infty}(K)}^2\int_{K \times K}|\phi(z)-\phi(w)|^2\, J(dzdw)\\
&\h + \|v\|_{L^{\infty}(K)}^2\int_{K \times K^c}|\phi(z)|^2\, J(dzdw).
\end{align*}
For the killing part, we have
\begin{align*}
\int_D |\phi v|^2\, dk \leq \| v\|_{L^{\infty}(K)}^2 \int_K \phi^2\, dk.
\end{align*}
Therefore, we obtain
\begin{align*}
\cE(\phi v, \phi v)\leq C\left( \mu_{\langle v \rangle}^c(K)+\int_{K \times K}|v(z)-v(w)|^2\, J(dzdw)+\| v \|_{L^{\infty}(K)}^2 \right)
\end{align*}
for some constant $C>0$ depending only on $\phi$, which is the desired inequality.
\end{proof}
For an open set $D$, the $\bdm|_{D}$-symmetric Hunt process associated with the part form $(\cE^D,\cF^D)$ is the part process $\bZ^D=(\Omega, \{Z_t^D\}_{t\geq 0}, \zeta^D, \{\Prob_z\}_{z\in D})$, where $\zeta^D=\tau_D$,
\begin{align*}
Z^D_t(\omega):=
\begin{cases}
Z_t(\omega),&\ t<\tau_D(\omega),\\
\partial,&\ t\geq \tau_D(\omega)
\end{cases}
\end{align*}
for $\omega \in \Omega$. We denote the set of smooth measures $\sS$, $\sS_0$, $\sS_{00}$ with respect to the part form $(\cE^D,\cF^D)$ by $\sS^D$, $\sS_0^D$, $\sS_{00}^D$, respectively. We also denote the $1$-potential of $\mu \in \sS_0^D$ by $U_1^D \mu$.
\begin{lem}\label{estmu}
Let $K$ be a compact set with $K\Subset D$ and $D_1$ a relatively compact open set with $K \Subset D_1 \Subset D$. Let $\nu_K$ be the $1$-equilibrium measure of $K$ with respect to the part form $(\cE^{D_1},\cF^{D_1})$.
For $v \in \cF^D$, we set
\begin{align*}
\cA_t:=\langle M^{[v],c}, M^{[v],c}\rangle_t+\int_0^t\int_K|\tilde v(Z_s)-\tilde v(w)|^2\, N(Z_s,dw)dH_s.
\end{align*}
Then it holds that
\begin{align*}
\mu_{\langle v \rangle}^c(K)+\int_{K \times K}|\tilde v(z)-\tilde v(w)|^2\, J(dzdw) \leq \Ex_{\nu_K} \left[ \int_0^{\tau_{D_1}}e^{-t}\, d \cA_t \right]
\end{align*}
\end{lem}
\begin{proof}
Since $\cA$ is a PCAF of $Z$, its stopped process $\cA^{\tau_{D_1}}$ is a PCAF of $Z^{D_1}$ and its corresponding Revuz measure $\mu_{\cA}$ is given by
\[
\mu_{\cA}(dz)=\mathbf{1}_{D_1}(z)\mu_{\langle v \rangle}^c(dz)+\mathbf{1}_{D_1}(z)\int_K|\tilde v(z)-\tilde v(w)|^2\, J(dzdw).
\]
Then by Theorem 2.3.15 of \cite{ChenFuku}, there exists a compact $\cE^{D_1}_1$-nest $\{F_k\}_{k=1}^{\infty}$ such that $\mathbf{1}_{F_k}\mu_{\cA} \in S_{00}^{D_1}$ and
\[
\mu_{\cA}\left(D_1\, \backslash \, \bigcup_{k=1}^{\infty}F_k\right)=0.
\]
For each $k$, by Lemma 4.1.5 of \cite{ChenFuku},
\[
\wilde{U^{D_1}_1(\mathbf{1}_{F_k}\mu_{\cA})}(z) = \Ex_z\left[ \int_0^{\tau_{D_1}}e^{-t}\mathbf{1}_{F_k}(Z_t)\, d\cA_t \right] \ \text{for q.e.}\ z\in D_1.
\]
Therefore,
\begin{align*}
\int_{D_1}\wilde{e_K}\mathbf{1}_{F_k} \, d\mu_{\cA} &= \cE^{D_1}_1(e_K, U^{D_1}_1(\mathbf{1}_{F_k}\mu_{\cA}))\\
&=\int_{D_1}\wilde{U^{D_1}_1(\mathbf{1}_{F_k}\mu_{\cA})}\, d\nu_K.
\end{align*}
Since $\wilde{e_K}=1$ q.e. on $K$, we obtain
\[
\mu_{\cA}(F_k \cap K)\leq \int_{D_1}\wilde{e_K}\mathbf{1}_{F_k} \, d\mu_{\cA} \leq \Ex_{\nu_K}\left[ \int_0^{\tau_{D_1}}e^{-t}\mathbf{1}_{F_k}(Z_t)\, d\cA_t \right].
\]
Letting $k \to \infty$, we obtain the desired estimate by the monotone convergence theorem.
\end{proof}

\begin{lem}\label{estdecom}
Let $M$ be a properly embedded compact Riemannian submanifold of $\mathbb{R}^d$. Let $u \in \cF^D_{loc}(M)$ be a quasi-harmonic map on $D$. Then for each relatively compact open subset $D_1 \Subset D$, the semimartingale $X=u(Z_{t\land \tau_{D_1}})$ has a special semimartingale decomposition
\[
X_t=X_0+L_t+C_t,
\]
where $L$ is a $\Prob_z$-square integrable martingale for q.e. $z\in E$ and $C$ is a continuous process of $\Prob_z$-integrable variation and $L$ and $C$ satisfy
\begin{align}
\Ex_{\nu}\left[ \int_0^{\tau_{D_1}} e^{-t}\, d[L,L]_t \right] &< \infty,\label{expbddL}\\
\Ex_{\nu}\left[ \int_0^{\tau_{D_1}}e^{-t}\, \left| dC_t \right| \right] &< \infty \label{expbddC}
\end{align}
for any $\nu \in S_{00}^{D_1}$.
\end{lem}
\begin{proof}
Let $\psi \in C_0(D)\cap \cF$ satisfy
\[
\psi = 1\ \text{on}\ \overline{D_1},\ 0\leq \psi \leq 1\ \text{on}\ E.
\]
Then $\psi u \in \cF^D$ and for $t\leq \tau_{D_1}$,
\begin{align*}
u(Z_t)-u(Z_0)&= \psi u(Z_t) - \psi u(Z_0) + \mathbf{1}_{\{ t\geq \tau_{D_1} \}}\left\{ \left( 1 - \psi \right)u(Z_{\tau_{D_1}}) - (1-\psi)u(Z_0) \right\}\\
&=M^{[\psi u]}_t + N^{[\psi u]}_t + \mathbf{1}_{\{ t\geq \tau_{D_1} \}}\left\{ \left( 1 - \psi \right)u(Z_{\tau_{D_1}}) - (1-\psi)u(Z_0) \right\}.
\end{align*}
Let
\begin{align*}
A_t&:=\mathbf{1}_{\{ t\geq \tau_{D_1} \}}\left\{ \left( 1 - \psi \right)u(Z_{\tau_{D_1}}) - (1-\psi)u(Z_0) \right\},\\
B_t&:=\int_0^{t\land \tau_{D_1}}\int_E \left\{ \left( 1 - \psi \right)u(w)-\left( 1 - \psi \right)u(Z_s)\right\}\, N(Z_s,dw)dH_s.
\end{align*}
Then $B_t$ is the $\Prob_z$-dual predictable projection of $A_t$ by the L\'evy system formula. We set
\begin{align*}
L_t&:=M^{[\psi u]}_{t\land \tau_{D_1}}+A_{t\land \tau_{D_1}}-B_{t\land \tau_{D_1}},\\
C_t&:=N^{[\psi u]}_{t\land \tau_{D_1}}+B_{t\land \tau_{D_1}}.
\end{align*}
Since $u(Z_{t\land \tau_{D_1}})$ is a $\Prob_z$-semimartingale, $N^{[\psi u]}_{t\land \tau_{D_1}}$ is a process of bounded variation.
Now it holds that
\begin{align*}
\Ex_{\nu}\left[ \int_0^{\tau_{D_1}} e^{-t}\, d[M^{[\psi u]},M^{[\psi u]}]_t \right]&= \Ex_{\nu}\left[ \int_0^{\infty} e^{-t} [M^{[\psi u]},M^{[\psi u]}]_{t\land \tau_{D_1}}\, dt \right]\\
&= \int_0^{\infty} \Ex_{\nu}\left[ \langle M^{[\psi u]},M^{[\psi u]}\rangle_{t\land \tau_{D_1}} \right]\, e^{-t}dt\\
&\leq \int_0^{\infty} (1+t)e^{-t}\, dt \| U_1^{D_1}\nu \|_{L^{\infty}(D_1)}\mu_{\langle \psi u\rangle}(D)\\
&=2\| U_1^{D_1}\nu \|_{L^{\infty}(D_1)}\mu_{\langle \psi u\rangle}(D),
\end{align*}
where we used Lemma 4.2.4 of \cite{ChenFuku}.
On the other hand, since $M$ is compact,
\begin{align*}
[A-B,A-B]_{\tau_{D_1}}=|\Delta A_{\tau_{D_1}}|^2\leq D_M^2 < \infty,
\end{align*}
where
\[
D_M:=\sup_{x\in M}|x|.
\]
Therefore,
\begin{align*}
\Ex_{\nu}\left[ \int_0^{\tau_{D_1}} e^{-t}\, d[A-B,A-B]_t \right]=\Ex_{\nu}\left[ e^{-\tau_{D_1}} |\Delta A_{\tau_{D_1}}|^2 \right]\leq D_M^2\nu(D_1)<\infty.
\end{align*}
Thus we get \eqref{expbddL}. In addition, since both $M^{[\psi u]}$ and $A-B$ are $\Prob_z$-square-integrable, $L$ is a $\Prob_z$-square-integrable martingale. By It\^o's formula for $\bi(X)$, $C_t$ coincides with
the process
\[
\frac{1}{2}\int_0^{t\land \tau_{D_1}} \partial_i \partial_j \bi(X_{s-})\, d[\bi^i(X),\bi^j(X)]^c_s + R^p_t,
\]
where $R^p$ is the $\Prob_z$-dual predictable projection of
\[
R_t:=\sum_{0<s \leq t\land \tau_{D_1}}\left\{ \bi(X_s)-\bi(X_{s-})-\partial_i \bi(X_{s-})\Delta \bi^i(X) \right\}
\]
since both of them are predictable bounded variation-parts of the special semimartingale decomposition of $X$. Since $M$ is compact,
\begin{align*}
\int_0^{t\land \tau_{D_1}} |\partial_i \partial_j \bi(X_{s-})|\, |d[\bi^i(X),\bi^j(X)]^c_s| &\leq c[X,X]^c_t,\\
\int_0^{t\land \tau_{D_1}}|dR_s| &\leq c \sum_{0<s \leq t\land \tau_{D_1}}|\Delta X_s|^2
\end{align*}
for some constant $c>0$ depending only on $M,\iota$.
In particular, this means that
\[
\Ex_z\left[\int_0^{t\land \tau_{D_1}}|dR^p_s| \right] \leq c \Ex_z\left[[X,X]_t\right]=c \Ex_z\left[[L,L]_t\right].
\]
Therefore, we get \eqref{expbddC} from \eqref{expbddL}.
\end{proof}
We note that since each $u_n$ is quasi-continuous on $D$, there exists an $\bdm$-polar set $\mathbf{N}$ such that
\begin{align}
|u_n(z)-u_m(z)|\leq \|u_n-u_m\|_{L^{\infty}(D_1)} \label{polarunif}
\end{align}
for all $z\in D_1 \backslash \mathbf{N}$ and $m,n \in \mathbb{N}$. Hereafter, we fix such $\mathbf{N}$.
\begin{lem}\label{excursion}
Under the assumptions of \tref{harmonicseq}, let $D_1 \Subset D$ be a relatively compact open set and $X^n:=u_n(Z)^{\tau_{D_1}}$.
Let $R>0$ be the radius obtained in \lref{Diffquad}.
Let $r, \ep>0$ be such that $r+2\ep <R$.
Let $\mathbf{N}$ be an $\bdm$-polar set satisfying \eqref{polarunif}.
Let $n_{\ep} \in \mathbb{N}$ be an integer such that for any $n \geq n_{\ep}$,
\[
\sup_{z\in \overline{D_1} \backslash \mathbf{N}}|u_n(z)-u_{n_{\ep}}(z)| \leq \ep.
\]
We denote $X=X^{n_{\ep}}$.
Define a sequence of stopping times $\{\sigma_j\}_{j=0}^{\infty}$ by
\begin{align*}
\sigma_0&:=0,\\
\sigma_{j+1}&:=\inf\{ t\geq \sigma_j \mid |X_t - X_{\sigma_j}| \geq r \}.
\end{align*}
Set
\[
\cN_t:=\sum_{j=0}^{\infty}\mathbf{1}_{\{\sigma_j < t\land \tau_{D_1}\}}.
\]
Then
\[
\Ex_{\nu}\left[ \int_0^{\infty}\cN_t e^{-t}\, dt \right]<\infty
\]
for any $\nu \in S_{00}^{D_1}$.
\end{lem}
\begin{proof}
Let $X=X_0+L+C$ be the decomposition obtained in \lref{estdecom}.
For each $j$,
\begin{align*}
r\mathbf{1}_{\{\sigma_{j+1} < t\land \tau_{D_1}\}} &\leq |X_{t\land \sigma_{j+1}}-X_{t\land \sigma_j}|\\
&\leq |L_{t\land \sigma_{j+1}}-L_{t\land \sigma_j}|+|C_{t\land \sigma_{j+1}}-C_{t\land \sigma_j}|.
\end{align*}
This means
\[
\frac{r}{2}\mathbf{1}_{\{\sigma_{j+1} < t\land \tau_{D_1}\}} \leq |L_{t\land \sigma_{j+1}}-L_{t\land \sigma_j}|\ \text{or}\ \frac{r}{2}\mathbf{1}_{\{\sigma_{j+1} < t\land \tau_{D_1}\}} \leq |C_{t\land \sigma_{j+1}}-C_{t\land \sigma_j}|.
\]
Hence, we have
\[
\mathbf{1}_{\{\sigma_{j+1} < t\land \tau_{D_1}\}} \leq \frac{4}{r^2}|L_{t\land \sigma_{j+1}}-L_{t\land \sigma_j}|^2+\frac{2}{r}|C_{t\land \sigma_{j+1}}-C_{t\land \sigma_j}|.
\]
Summing up over $j=0,1,\dots$, we have
\[
\cN_t-1 \leq \frac{4}{r^2}\sum_{j=0}^{\infty}|L_{t\land \sigma_{j+1}}-L_{t\land \sigma_j}|^2 + \frac{2}{r}\sum_{j=0}^{\infty}|C_{t\land \sigma_{j+1}}-C_{t\land \sigma_j}|.
\]
By the orthogonality of the square-integrable martingale $L$,
\[
\sum_{j=0}^{\infty}\Ex_z\left[ |L_{t\land \sigma_{j+1}}-L_{t\land \sigma_j}|^2 \right] \leq \Ex_z\left[ [L,L]_t \right].
\]
We also have
\[
\sum_{j=0}^{\infty}\Ex_z\left[|C_{t\land \sigma_{j+1}}-C_{t\land \sigma_j}|\right] \leq \Ex_z \left[ \int_0^{t\land \tau_{D_1}}|dC_t| \right].
\]
Therefore, we have
\[
\Ex_{\nu}\left[ \cN_t \right] \leq \nu(D_1) + \frac{4}{r^2} \Ex_{\nu}\left[ [L,L]_t \right] + \frac{2}{r} \Ex_{\nu} \left[ \int_0^{t\land \tau_{D_1}}|dC_t| \right].
\]
Integrating over $t\in [0,\infty)$ with respect to the exponential distribution, we obtain
\[
\Ex_{\nu}\left[ \int_0^{\infty}\cN_t e^{-t}\, dt \right] \leq \nu(D_1) +  \frac{4}{r^2} \Ex_{\nu}\left[\int_0^{\tau_{D_1}} e^{-t}\, d[L,L]_t \right] + \frac{2}{r} \Ex_{\nu} \left[ \int_0^{\tau_{D_1}}e^{-t}\, |dC_t| \right]<\infty
\]
by \lref{estdecom}.
\end{proof}
\begin{proof}[Proof of \tref{harmonicseq}]
Let $\phi \in C_0(D)\cap \cF$ and set $K=\mathrm{supp}[\phi]$. Take a relatively compact open subset $D_1$ with $K \Subset D_1 \Subset D$. Since each $u_n$ is bounded, there exists $u_n' \in \cF^D(\mathbb{R}^d)\cap L^{\infty}(E)$ such that
\[
u_n'=u_n\ \bdm \text{-a.e. on}\ D_1.
\]
Then $\phi u_n=\phi u_n'$ $\bdm$-a.e. on $E$. Hence by \lref{estE}, there exists a constant $C>0$ depending only on $\phi$ such that
\begin{align}
\cE(\phi(u_n-u_m), \phi(u_n-u_m)) &= \cE(\phi v_{m,n}, \phi v_{m,n})\nonumber \\
&\leq C \left(\mu_{\langle v_{m,n} \rangle}^c(K) + \int_{K \times K}|v_{m,n}(z)-v_{m,n}(w)|^2\, J(dzdw)\right. \nonumber \\
&\h \left.+ \| v_{m,n} \|_{L^{\infty}(K)}^2 \right),\label{Ephiumn}
\end{align}
where we set $v_{m,n}=u_n'-u_m'$.
By \lref{estmu},
\begin{align}
\mu_{\langle v_{m,n} \rangle}^c(K) + \int_{K \times K}|v_{m,n}(z)-v_{m,n}(w)|^2\, J(dzdw) &\leq \Ex_{\nu_K}\left[ \int_0^{\tau_{D_1}}e^{-t}\, d \cA_t \right]\nonumber \\
&=\int_0^{\infty} \Ex_{\nu_K}\left[ \cA_{t\land \tau_{D_1}} \right] e^{-t}\, dt\nonumber \\
&\leq \int_0^{\infty} \Ex_{\nu_K}\left[ Q^{m,n}_t \right] e^{-t}\, dt,\label{bddmuK}
\end{align}
where $X^n:=u_n(Z)^{\tau_{D_1}}$ and $Q^{m,n}_t:=[ X^n-X^m, X^n-X^m]_{t}$
since it holds that
\begin{align*}
\Ex_z\left[ \cA_{t\land \tau_{D_1}} \right] &= \Ex_z\left[ \langle M^{[v_{m,n}],c},M^{[v_{m,n}],c} \rangle_{t\land \tau_{D_1}} \right] + \Ex_z \left[ \sum_{0<s\leq t\land \tau_{D_1}}|\Delta \tilde v_{m,n}(Z_s)|^2 \mathbf{1}_K(Z_s) \right]\\
&\leq \Ex_z\left[ Q^{m,n}_t \right]
\end{align*}
for q.e. $z\in D_1$.
Let $R, r, \ep>0$ be as in \lref{excursion}.
Let $\mathbf{N}$ be an $\bdm$-polar set such that for all $m,n\in \mathbb{N}$,
\begin{align*}
\begin{cases}
&u_n'=u_n\ \text{on}\ D_1 \backslash \mathbf{N},\\
&|u_n(z)-u_m(z)| \leq \| u_n-u_m\|_{L^{\infty}(D_1)}\ \text{for}\ z\in D_1 \backslash \mathbf{N},\\
&u_n \to u\ \text{pointwise on}\ E \backslash \mathbf{N},\\
&u_n(Z)^{\tau_{D_1}}\ \text{is a}\ \Prob_z \text{-martingale on}\ M\ \text{for all}\ z\in E\backslash \mathbf{N}.
\end{cases}
\end{align*}
Let $n_{\ep} \in \mathbb{N}$ be a positive integer such that
\[
\|u_n-u_{n_{\ep}}\|_{L^{\infty}(\overline{D_1})} \leq \ep
\]
for all $n \geq n_{\ep}$.
Let $\{\sigma_j\}_{j=0}^{\infty}$ be a sequence of stopping times given in \lref{excursion}. Set $\Lambda_j:=\{ \sigma_j < t\land \tau_{D_1} \}$ for fixed $t\geq 0$. Then for each $n$, the process $\{ X^n_{(\sigma_j+s)\land (t\land \tau_{D_1} \land \sigma_{j+1})} \}_{s\geq 0}$ is an $M$-valued $(\Prob_z(\cdot \mid \Lambda_j), \{\cF_{\sigma_j+s}\}_{s\geq 0})$-martingale with an end point and satisfies
\begin{align*}
\sup_{\sigma_j \leq s < \sigma_{j+1}\land t \land \tau_{D_1}}|X^n_s-X^n_{\sigma_j}| &\leq \sup_{\sigma_j \leq s < \sigma_{j+1}\land t \land \tau_{D_1}}\left(|X^n_s-X^{n_{\ep}}_s| + |X^{n_{\ep}}_s-X^{n_{\ep}}_{\sigma_j}|  \right) + |X^{n_{\ep}}_{\sigma_j}-X^n_{\sigma_j}| \\
&\leq 2 \ep + r \\
&<R
\end{align*}
for all $n \geq n_{\ep}$.
Therefore, for $m,n \geq n_{\ep}$, by \lref{Diffquad},
\begin{align*}
\Ex_z \left[ Q^{m,n}_{t\land \sigma_{j+1}}- Q^{m,n}_{t\land \sigma_j};\, \Lambda_j \right]&=\Prob_z(\Lambda_j)\cdot \Ex_z \left[ Q^{m,n}_{t\land \sigma_{j+1}}- Q^{m,n}_{t\land \sigma_j} \mid \Lambda_j \right]\\
&\leq C_1 \Prob_z(\Lambda_j)\cdot \Ex_z \left[ \left. \left( \sup_{s \in [t\land \sigma_j, t\land \sigma_{j+1}]}\left|X^m_s - X^n_s \right| \right)^2 \right| \Lambda_j \right]^{\frac{1}{{2}}}\\
&= C_1 \Prob_z(\Lambda_j)^{\frac{1}{2}}\cdot \Ex_z \left[ \left( \sup_{s \in [t\land \sigma_j, t\land \sigma_{j+1}]}\left|X^m_s - X^n_s \right| \right)^2 ;\, \Lambda_j \right]^{\frac{1}{{2}}}
\end{align*}
for some constant $C_1>0$ provided in \lref{Diffquad}.
The above inequality holds even if $\Prob_z(\Lambda_j)=0$.
Integrating both sides with respect to $\nu_K$, we get
\begin{align*}
\Ex_{\nu_K} &\left[ Q^{m,n}_{t\land \sigma_{j+1}}- Q^{m,n}_{t\land \sigma_j};\, \Lambda_j \right]\\
&\leq C_1 \int_{D_1} \Prob_z(\Lambda_j)^{\frac{1}{2}}\cdot \Ex_z \left[ \left( \sup_{s \in [t\land \sigma_j, t\land \sigma_{j+1}]}\left|X^m_s - X^n_s \right| \right)^2 ;\, \Lambda_j \right]^{\frac{1}{{2}}}\, \nu_K(dz)\\
&\leq C_1 \Prob_{\nu_K}(\Lambda_j)^{\frac{1}{2}}\cdot \Ex_{\nu_K} \left[ \left( \sup_{s \in [t\land \sigma_j, t\land \sigma_{j+1}]}\left|X^m_s - X^n_s \right| \right)^2 ;\, \Lambda_j \right]^{\frac{1}{2}}
\end{align*}
by the Cauchy-Schwarz inequality. Taking the sum over $j$, we get
\begin{align*}
\Ex_{\nu_K} \left[ Q^{m,n}_t \right] &\leq C_1 \left(\sum_{j=0}^{\infty}\Prob_{\nu_K}(\Lambda_j) \right)^{\frac{1}{2}} \cdot \left( \sum_{j=0}^{\infty} \Ex_{\nu_K} \left[ \left( \sup_{s \in [t\land \sigma_j, t\land \sigma_{j+1}]}\left|X^m_s - X^n_s \right| \right)^2 ;\, \Lambda_j \right] \right)^{\frac{1}{2}}\\
&\leq C_1 \Ex_{\nu_K}[\cN_t]^{\frac{1}{2}}\cdot \Ex_{\nu_K} \left[ \cN_t \left( \sup_{s \in [0, t]}\left|X^m_s - X^n_s \right| \right)^2 \right]^{\frac{1}{2}}
\end{align*}
by again the Cauchy-Schwarz inequality.
By further integrating both sides with respect to the exponential distribution, we obtain
\begin{align*}
\int_0^{\infty}\Ex_{\nu_K} \left[Q^{m,n}_t \right]e^{-t}\, dt &\leq C_1 \left( \int_0^{\infty} \Ex_{\nu_K}[\cN_t] e^{-t}\, dt\right)^{\frac{1}{2}} \cdot\\
&\h \left( \int_0^{\infty} \Ex_{\nu_K} \left[ \cN_t \left( \sup_{s \in [0, t]}\left|X^m_s - X^n_s \right| \right)^2 \right] e^{-t}\, dt \right)^{\frac{1}{2}}.
\end{align*}
By \lref{excursion},
\[
\int_0^{\infty} \Ex_{\nu_K}[\cN_t] e^{-t}\, dt<\infty.
\]
On the other hand, since
\[
u_m-u_n\to 0\ \text{in}\ L^{\infty}(\overline{D_1})\ \text{and q.e. on}\ E,
\]
we have
\[
\sup_{s \in [0, l\land \tau_{D_1}]}\left|X^m_s - X^n_s \right| \to 0\ \Prob_z\text{-a.s. for q.e.}\ z\in E \ \text{and all}\ l \in \mathbb{N}.
\]
By the dominated convergence theorem, it holds that
\[
\lim_{m,n \to \infty} \int_0^{\infty} \Ex_{\nu_K} \left[ \cN_t \left( \sup_{s \in [0, t]}\left|X^m_s - X^n_s \right| \right)^2 \right] e^{-t}\, dt =0.
\]
Hence, by \eqref{bddmuK} and \eqref{Ephiumn}, $\{\phi u_n\}_{n \in \mathbb{N}}$ is an $\cE_1$-Cauchy sequence and $\phi u_n$ converges to $\phi u$ in $\cE_1$. In particular, $u \in \cF_{loc}^D(\mathbb{R}^d)$ and we can take a modification of $u$ which is Borel measurable on $E$ and quasi-continuous on $D$. Since $u_n \to u$ q.e. on $E$ by assumption and $M$ is closed in $\mathbb{R}^d$, $u(z)\in M$ q.e.
Moreover, since $u_n \to u$ in $L^{\infty}_{loc}(D)$ and pointwise on $E$ except for an $\bdm$-polar set, it holds that
\[
\sup_{s\in [0,t\land \tau_{D_1}]}|u(Z_s)-u_n(Z_s)| \leq \|u-u_n\|_{L^{\infty}(D_1)}+\mathbf{1}_{\{\tau_{D_1}\leq t\}}|u(Z_{\tau_{D_1}})-u_n(Z_{\tau_{D_1}})| \to 0
\]
as $n \to \infty$ $\Prob_z$-a.s. for q.e. $z\in E$. Therefore, by \cref{equitopo}, $u(Z)^{\tau_{D_1}}$ is an $M$-valued $(\Prob_z,\{\cF^Z_t\}_{t\geq 0})$-martingale with the killing time $\hat{\zeta}$ and the end point $0$ for q.e. $z\in E$. Hence, $u$ is a quasi-harmonic map.
\end{proof}


\section*{Statements and Declarations}
{\bf Competing interests:} No conflict is related to this article, and the author has no relevant financial or non-financial interests to disclose.\\
{\bf Ethical Approval:} Not applicable to this article.\\
{\bf Funding:} This work was supported by JSPS KAKENHI Grant Number 24K22827, 25K17264.\\
{\bf Data Availability Statement:} Data sharing is not applicable to this article as no datasets were generated or analyzed during the current study.

\end{document}